\documentclass[12pt]{article}
\usepackage[utf8]{inputenc}
\usepackage{amsmath, amssymb, amsthm, enumerate}
\usepackage{graphicx, xcolor}
\usepackage{mathtools}
\usepackage{caption}

\newtheorem{theorem}{Theorem}[section]
\newtheorem{lemma}[theorem]{Lemma}

\newtheorem{cor}[theorem]{Corollary}

\def\ff{{\cal F}}
\def\rr{{\mathbb R}}
\def\zz{{\mathbb Z}}
\def\sik{{\rr}^2}
\def\su{\subset}
\def\al{\alpha}
\def\be{\beta}
\def\ga{\gamma}
\def\de{\delta}
\def\ep{\varepsilon}
\def\si{\sigma}
\def\Si{\Sigma}
\def\la{\lambda}

\def\stb{,\ldots ,}

\def\noi{\noindent}

\def\ol{\overline}
\def\proof{\noi {\bf Proof.}\ }
\def\szog{\measuredangle}

\title{Quadrilateral reptiles}

\author{Mikl\'os Laczkovich\thanks{Partially supported by Hungarian Scientific
Foundation grant no. K124749.}}

\begin{document}

\maketitle 

\begin{abstract}
A polygon $P$ is called a reptile, if it can be decomposed into $k\ge 2$
nonoverlapping and congruent polygons similar to $P$.
We prove that if a cyclic quadrilateral is a reptile, then it is a trapezoid.
Comparing with results of U. Betke and I. Osburg we find that every convex
reptile is a triangle or a trapezoid.
  
\end{abstract}

\insert\footins{\footnotesize{MSC code: 52C20, 52A10}}
\insert\footins{\footnotesize{Key words: tilings, reptiles, quadrilaterals}}

\section{Introduction and main results}

A polygon $P$ is called a {\it $k$-reptile}, if it can be decomposed into $k$
nonoverlapping and congruent polygons similar to $P$. For example, every
triangle and every parallelogram is a $4$-reptile. A polygon $P$ is a {\it
reptile}, if it is a $k$-reptile for some $k\ge 2$.

C. D. Langford presented several reptiles, including three $4$-reptile
trapezoids, and raised the problem of classifying all reptiles \cite{L}.
It was shown by G. Valette and T. Zamfirescu in \cite{VZ} that every convex
$4$-reptile is a triangle, a parallelogram or one of Langford's trapezoids.

Non-convex reptiles are abundant, and their characterization seems to be
hopeless (see \cite[Section 3.2]{O}). As for convex reptiles, it is known
that they must be triangles or quadrilaterals (see \cite{B} and
\cite[Satz 2.23]{O}). A large step towards the classification of convex
reptiles was made by I. Osburg. She proved in \cite[Satz 2.9 and
Folgerung 3.2]{O} that every quadrilateral reptile (convex or not) is either
a trapezoid or a cyclic quadrilateral. In this note our aim is to prove
the following.

\begin{theorem}\label{t1}
If a cyclic quadrilateral is a reptile, then it is a trapezoid.
\end{theorem}

Comparing with Osburg's results we obtain the following corollaries
(the first one was conjectured by J. Doyen and M. Landuyt in \cite{DL}).
\begin{cor}\label{c1}
Every convex reptile is a triangle or a trapezoid.
\end{cor}

\begin{cor}\label{c2}
Every quadrilateral reptile (convex or not) is a trapezoid.
\end{cor}  
The characterization of reptile trapezoids remains open.

It follows from Corollary \ref{c1} that {\it the class of self-similar
quadrilaterals is strictly larger than the class of quadrilateral reptiles.}
Indeed, one can show that every right kite (a kite with two right
angles) is self-similar. On the other hand, according to Corollary \ref{c1}
a kite is not a reptile (!), unless it is a rhombus. Therefore, a right kite
is self-similar but is not a reptile, unless it is a square.

The characterization of self-similar quadrilaterals
seems to be open as well. It is also unsolved if self-similar convex pentagons
exist or not. More is known about the even larger class of self-affine polygons.
It is proved in \cite{HR} that every convex quadrilateral is self-affine, and
that self-affine convex pentagons do exist.

The structure of our proof of Theorem \ref{t1} is the following. Let $Q$ be a
cyclic quadrilateral, and suppose that $Q$ is a reptile. Then there is a polygon
$Q'$ similar to $Q$ such that $Q'$ can be tiled
with congruent copies of $Q$. In the next section we show that already the
existence of a tiling in a neighbourhood of the smallest angle of $Q'$
implies that either $Q$ is a trapezoid, or $Q$ belongs to one of three
families, each one depending on a single continuous parameter (see Theorem
\ref{t2}).

In Section 3 we show that the quadrilaterals belonging to the
first two families are not reptiles unless they are trapezoids. The proofs are
relatively easy, using either number theoretical conditions (similar to the
argument of \cite{SWW}) or a local argument similar to the proof of Theorem
\ref{t2}.

The third family, denoted by $\ff _3$, consists of the quadrilaterals with
sides $a,b,c,d$ such that $c=d$ and the angles included between the sides $a,b$
and between $c,d$ are right angles. Assuming that $c=d=1$, the elements of
$\ff _3$ are determined by the length of $a$. The proof that the elements of
$\ff _3$ are not reptiles (apart from the square) occupies the last
two sections.

\begin{figure}
\includegraphics[width=3in]{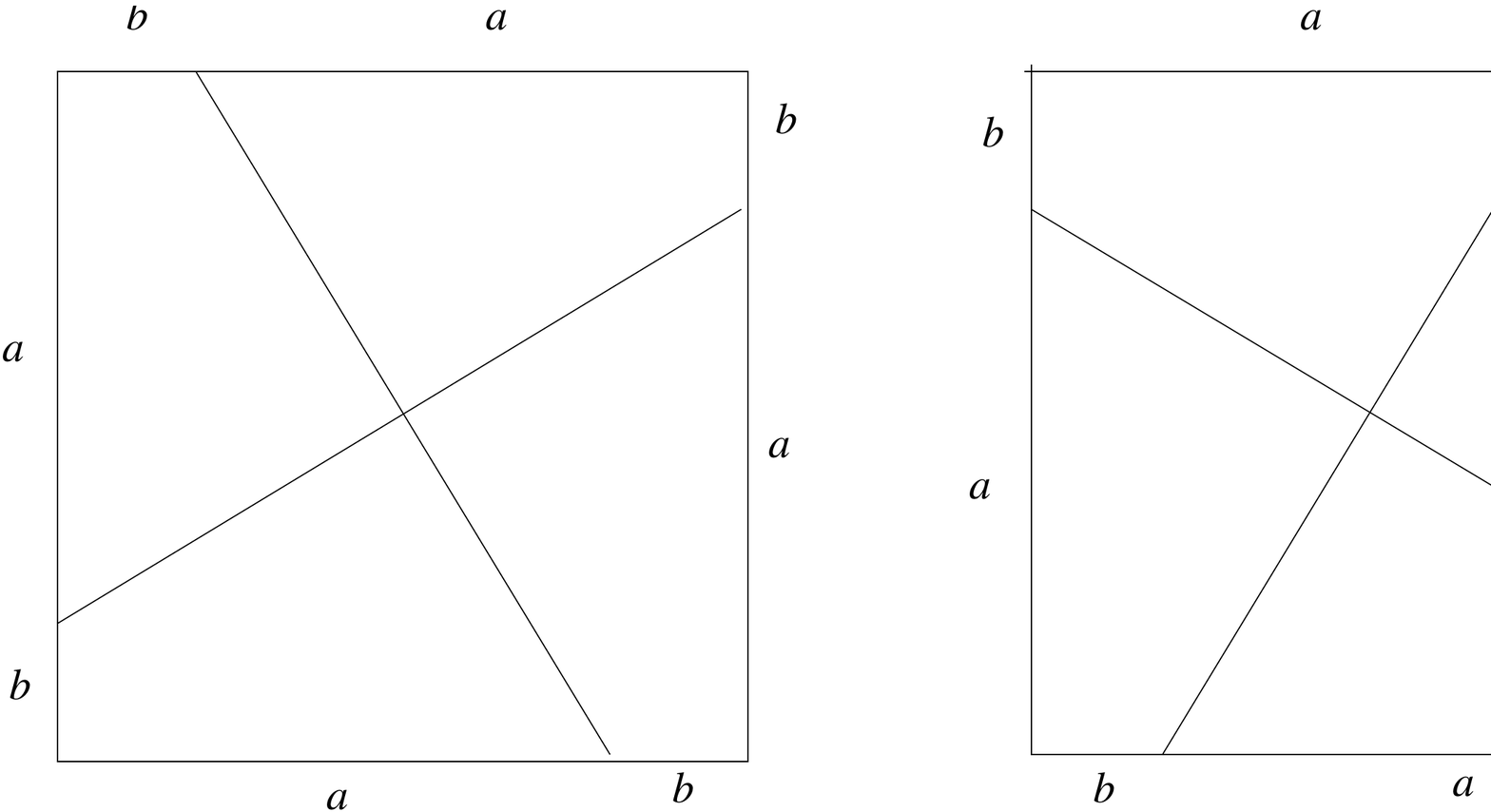}
\caption{}
\label{f0}
\end{figure}

There are two reasons why this proof is more difficult than ruling out the
other two families. First, the quadrilaterals belonging to $\ff _3$
tile the square, as Figure \ref{f0} shows. 
Consequently, they tile the whole plane in such a way that several local tilings
can be extended indefinitely, so the argument used in Theorem \ref{t2} doesn't
work. Also, the number theoretical arguments fail in the case when the side
lengths of the quadrilateral are rational. There are infinitely many such
quadrilaterals, corresponding to the rational solutions of $a^2 +b^2 =2$.
The special case $a=7/5$, $b=1/5$ (when $a=1+2b$) is particularly difficult to
handle.

In order to deal with the family $\ff _3$ we need a global view of the
possible tilings.
This is provided by Theorem \ref{t3}. It says that if a quadrant is tiled with
congruent copies of a quadrilateral belonging to $\ff _3$, then the tiling
is obtained by first tiling the quadrant by squares of side length $a+b$, and
then tiling these squares as in Figure \ref{f0}. From Theorem \ref{t3} it is
easy to infer that the elements of $\ff _3$ are not reptiles (except the
square).

We conclude this section with some preliminary remarks. The segment with
endpoints $A$ and $B$ is denoted by $AB$. We denote the length of the
segment $AB$ by $\ol{AB}$.

Let $Q$ be a cyclic quadrilateral, and suppose that $Q$ is not a trapezoid.
Let $\al ,\be ,\ga ,\de$ be the angles of $Q$ listed counterclockwise and
such that $\al$ is (one of) the smallest of the angles. Since $\de =\pi -\be$
and $\ga =\pi -\al$, this implies that $\ga$ is (one of) the largest of the
angles; that is, $\al \le \be \le \ga$ and $\al \le \de \le \ga$.

Since  $Q$ is not a trapezoid, we have $\al \ne \be$ and $\al \ne \de$.
Therefore, we have $\al <\be <\ga$ and $\al <\de <\ga$. 

In the sequel we fix the labeling of the vertices, angles and sides of $Q$ as
follows. The vertices $A,B,C,D$ of $Q$ are listed such that the angles of $Q$
at the vertices $A,B,C,D$ are $\al ,\be ,\ga ,\de$. The lengths of the sides of $Q$ are $a,b,c,d$ such that $\ol{AB}=a$, $\ol{BC}=b$, $\ol{CD}=c$ and
$\ol{DA}=d$. We also assume that $a\ge d$. (The case $d\ge a$ can be
reduced to $a\ge d$ if we turn to a reflected copy of $Q$ instead of $Q$,
and swap $\be$ and $\de$.) 

\begin{lemma}\label{l1}
We have $a\ge d>b$ and $a>c$.
\end{lemma}
\proof Since $\al < \ga =\pi -\al$, we have $\al < \pi /2< \ga$.
If $d\ge \ol{BD}$ then, as angle of the triangle $BCD$ at the vertex $C$ equals
$\ga > \pi /2$, we have $b<\ol{BD}\le d$.

Next suppose $d<\ol{BD}$. Let $K$ denote the circumscribed circle
of $Q$. Let $ABD\szog =\ep$. Since $d<\ol{BD}$, we have $\ep <\al$. On the
other hand we have $\be > \al$, and thus there is a point $E$ in the subarc
of $K$ with end points $C$ and $D$ such that $ABE\szog =\al$. Then $ABED$ is a
trapezoid, and thus $\ol{BE}=d$. Now the angle of the triangle $BCE$ at the
vertex $C$ is greater than $\ga > \pi /2$, and thus $b<\ol{BE}=d$.

The proof of $a>c$ is similar. (Note that the proof of $d>b$ did not use
$a\ge d$.) \hfill $\square$

\begin{lemma}\label{l2}
If $c\ge d$, then $\be <2\al$.
\end{lemma}
\proof Let $ABD\szog =\ep$ and $DBC\szog =\zeta$. Since $Q$ is cyclic, we have
$DAC\szog =\zeta$. Also, $c\ge d$ implies $\ep \le \zeta$, and thus
$\be =\ep +\zeta \le 2\zeta <2\al$. \hfill $\square$

\section{Reduction to one dimensional families}
Let $Q$ be a cyclic quadrilateral, and suppose that $Q$ is not a trapezoid.
We use the labeling of the angles, vertices and sides of $Q$ as described in
the previous section. In this section we prove the following.
\begin{theorem}\label{t2}
If $Q$ is a reptile, then one of the following statements is true.
\begin{enumerate}[{\rm (i)}]
\item $\al =\pi /3$, $\be =\de =\pi /2$ and $\ga =2\pi /3$.
\item $b=c$, $\al =\pi /3$ and $\ga =2\pi /3$.
\item $c=d$, $\be =\pi/2$ and $\de =\pi/2$.
\end{enumerate}
\end{theorem}
\proof
Since $Q$ is a reptile, there is a quadrilateral $Q'$ similar to $Q$ such
that $Q'$ can be tiled with at least two congruent copies of $Q$.
We fix such a tiling, and fix the labeling of the angles and sides of the
tiles. If $T$ is a tile with vertices $X,Y,V,W$, then $|XY|, |YV|, |VW| , |WX|$
denote the labels of the sides of $T$. That is, if $|XY|=a$ then $\ol{XY}=a$
and there are two possibilities: either $|YV|=b, |VW|=c , |WX|=d$ and the
angles of $T$ at $X,Y,V,W$ are $\al ,\be , \ga ,\de$, or $|YV|=d, |VW|=c ,
|WX|=b$ and the angles of $T$ at $X,Y,V,W$ are $\al ,\de , \ga ,\be$. Note
that if $|XY|=a$ then $\ol{XY}=a$, but the converse is not necessarily true;
it can happen, e.g., that $\ol{XY}=a$ but $|XY|=d$. In this case, however, we
must have $a=d$.

Labeling the angles has similar consequences. When we say, e.g., that a tile $T$
has angle $\be$ at the vertex $Y$, then it implies that either $|XY|=a$
and $|YV|=b$, or $|XY|=b$ and $|YV|=a$. 

The angle $\al$ of $Q'$ is packed with one single tile. We may assume that
$Q$ itself is this tile; that is, $A$ is a vertex of $Q'$, and the sides
$AB$ and $AD$ lie on consecutive sides of $Q'$.

Let $U$ denote the union of the boundaries of the tiles. We say that a segment
$S\su U$ is maximal if, for every segment $S'$ with $S\su S' \su U$
we have $S'=S$. We consider three cases and several subcases.

\noi
{\bf Case I:} {\it the segment $DC$ is not maximal.} Then there is a segment
$DE \su U$ such that $C$ is an inner point of $DE$, and there is a tile $T$
having $C$ as a vertex, having angle $\pi -\ga =\al$ at $C$, and having a
side $CF$ on the segment $CB$. However, since $T$ and $Q$ are congruent,
$\ol{CF}$ equals either $a$ or $d$, and thus $d\le \ol{CB}=b$, which contradicts
Lemma \ref{l1}.

\noi
{\bf Case II:} {\it the segment $DC$ is maximal, and the segment $BC$ is also
maximal.} Then $C$ is a common vertex of the tiles $Q_0 ,Q_1 \stb Q_k$ such that
$Q_0 =Q$, and the sum of the angles $Q_i$ at $C$ equals $2\pi$. 

Let $Q_1$ be the tile with a vertex at $C$ and having a side $CE\su CB$.
Since $a\ge d>b$ by Lemma \ref{l1}, we have $|CE|=b$ or $c$. 

\noi
{\bf Case II.1:} $|CE|=c$ and $c<b$. In this case
the segment $BC$ is packed with $k\ge 2$ segments of length $c$, and each of
them is a side of a tile. That is, there is a partition $C=E_0 ,E=E_1 \stb
E_k =B$ and there are tiles $T_1 =Q_1 , T_2 \stb T_k$ such that $E_{i-1} E_i$ is
a side of $T_i$ labeled as $c$. The angle of $T_i$ at the vertices $E_{i-1}$
and $E_i$ equals either $\de$ or $\ga$. Since $\de +\ga >\al +\ga =\pi$,
$T_i$ cannot have angle $\ga$ at $E_{i-1}$ if $i\ge 2$, and cannot have angle
$\ga$ at $E_{i}$ if $i<k$. This is possible only if $k=2$, $T_1 =Q_1$ has
angle $\ga$ at $C=E_0$ and $T_2$ has angle $\ga$ at $B=E_2$. In this case,
however, the sum of the angles of $Q_0$ and $T_2$ at $B$ equals $\be +\ga >\pi$,
which is impossible.

\noi
{\bf Case II.2:} {\it $|CE|=b$, or $|CE|=c$ and $c=b$.} Then
$E=B$, and the angle of $Q_1$ at $B$ is $\be$, $\de$ or $\ga$. Since the angle
of $Q_0$ at $B$ equals $\be$ and $\be +\ga >\al +\ga =\pi$, the angle of $Q_1$
at $B$ is $\be$ or $\de$. Therefore, the angle of $Q_1$ at $C$ is $\ga$.
Note that the angle of $Q_1$ at $B$ can be $\de$ only if $c=b$.

The sum of the angles $Q_i$ ($i=2\stb k$) at $C$ equals $2\pi -2\ga
=2\al$. Since $\al$ is the smallest angle, this implies $k=2$ or $k=3$.
In the latter case the angles of $Q_2$ and $Q_3$ at $C$ equal $\al$. This,
however, contradicts the condition that $BC$ is maximal.

Thus $k=2$, and we have $2\ga +\eta =2\pi$, where $\eta$ is the angle of
$Q_2$ at $C$. Then $\eta =2\al$, where $\eta$ equals one of $\be , \ga ,\de$. 
Let $FC$ be the side of $Q_2$ lying on the segment $DC$.
Since $\ol{FC}\le c<a$, we have $|FC|=b$, $c$ or $d$.

\noi
{\bf Case II.2.1:} $\eta =\be =2\al$. Then the only possible label for
$FC$ is $b$. By Lemma \ref{l2}, $\be =2\al$ implies $c<d$.
Now $|FC|=b$ implies that the angle of $Q_2$ at $F$ is $\ga$. Since $\ga >\be$,
it follows that $F\ne D$, and there is a tile $Q_3$ having a vertex at $F$
and angle $\al$ at $F$. If $GF$ is the side of $Q_3$ lying on the segment $DF$,
then $\ol{GF}=a$ or $d$. However, we have $\ol{DF}<c<d\le a$, a contradiction.

\noi
{\bf Case II.2.2:} $\eta =\ga$. Then we have $3\ga =2\pi$, $\ga =2\pi /3$ and
$\al =\pi /3$. Now $2\be +\al >3\al =\pi$ and $\be +\de +\al >3\al =\pi$, and
thus the sum of the angles of $Q_0$ and $Q_1$ at $B$ equals $\pi$.
If the angle of $Q_1$ at $B$ equals $\be$, then $\be =\pi /2$, $\de =\pi /2$,
and we have case (i) of the lemma.

If the angle of $Q_1$ at $B$ equals $\de$, then necessarily $c=b$, 
and we have case (ii) of the lemma.

\noi
{\bf Case II.2.3:} $\eta =\de =2\al$. Then we have $|FC|=c$ or $d$.
By $\al <\be =\pi -\de =\pi -2\al$ we obtain $\al <\pi /3$.

\begin{figure}
\includegraphics[width=3in]{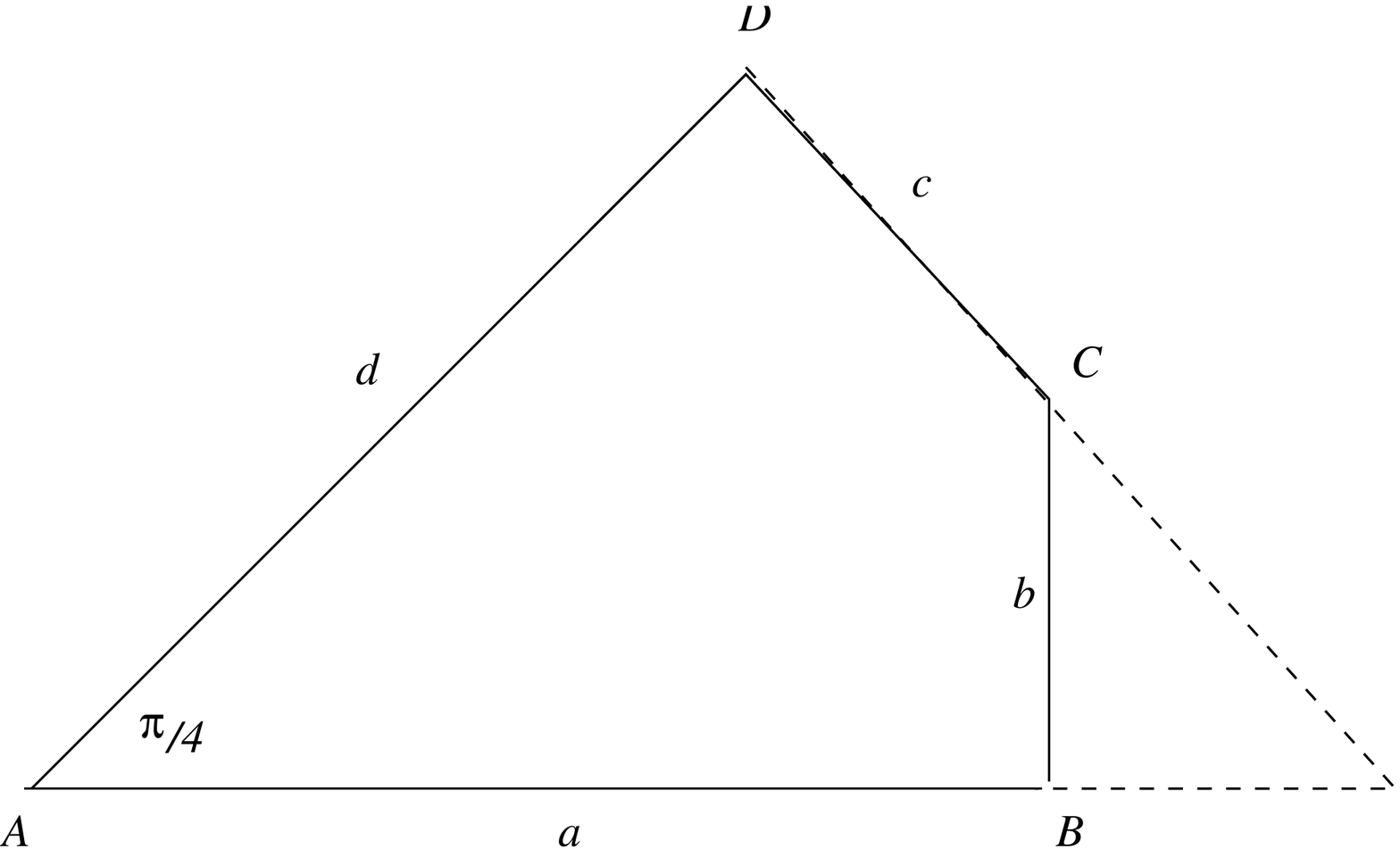}
\caption{}
\label{f1}
\end{figure}

\noi
{\bf Case II.2.3.1:} $c\ge d$. Then $c>b$, and thus the angle of $Q_1$ at
$B$ equals $\be$. We have $2\be +\al =2\pi -4\al +\al =2\pi -3\al >\pi$, and
thus the sum of the angles of $Q_0$ and $Q_1$ at $B$ equals $\pi$. Thus
$2\be =\pi$, $\be =\de =\pi /2$, $\al =\pi /4$. In this case, however, we have
$d>c$, a contradiction (see Figure \ref{f1}).

\noi
{\bf Case II.2.3.2:} $c<d$. Then $|FC|=c$, $F=D$, and the angle of $Q_2$
at $D$ must be $\ga$. However, the sum of the two angles
at $D$ equals $\de +\ga >\de +\be =\pi$, which is impossible.

\noi
{\bf Case III:} {\it the segment $DC$ is maximal, and the segment $BC$ is not
maximal.} Then $C$ is a common vertex of the tiles $Q$ and $R$ such that
$R$ has angle $\al$ at $C$, and there is a point $F$ in the segment
$DC$ such that $FC$ is a side of $R$. Thus $|FC|=a$ or $d$. Since
$a>c$ by Lemma \ref{l1}, $|FC|=a$ is impossible.
Thus $|FC|=d$ and the angle of $R$ at $F$ equals $\de$. We also have
\begin{equation}\label{e2}
c\ge d \ \text{and} \ \be <2\al .
\end{equation}
The second inequality follows from Lemma \ref{l2}.

\noi
{\bf Case III.1:} $F=D$. Then $c=d$ and the angle of $R$ at $D$ equals $\de$.
Now the angle of $Q$ at $D$ is also $\de$, hence $2\de \le \pi$. If
$2\de =\pi$, then $\de =\be =\pi/2$ and we get (iii) of the lemma.

If $2\de <\pi$, then there are other tiles at $D$. However, by \eqref{e2}
we have $2\de +\al >\de +2\al >\de +\be =\pi$, which is impossible.

\noi
{\bf Case III.2:} $F\ne D$. Then there is a partition $D=F_0 \stb F_n =C$
of the segment $DC$ and there are tiles $R_1 \stb R_n$ such that $n\ge 2$,
$F_{i-1} F_i$ is a side of $R_i$ for every $i=1\stb n$, $F_{n-1}=F$ and $R_n =R$.
It is clear that $|F_{i-1} F_i |$ is $b$ or $d$ for every $i$.
Since $|F_{n-1}C|=d$, we have $d<c$. Comparing with Lemma \ref{l1} we obtain
\begin{equation}\label{e3}
a>c>d>b.
\end{equation}
We denote by $\la _i$ and $\mu _i$ the angle of $R_i$ at $F_{i-1}$ and at
$F_i$, respectively.

We prove that $\be =\de =\pi /2$. Suppose this is not true.
Since $|F_{n-1} F_n |=|FC|=d$ and $\mu _n =\al$, it follows that
$\la _n =\de$. Thus $\mu _{n-1} \ne \ga$ by $\de +\ga >\pi$. Suppose
$\mu _{n-1} = \al$. Since $\al +\de <\pi$, there must be a third tile having
a vertex at $F_{n-1}$. However, $2\al +\de >\de +\be =\pi$, so this case is
impossible. Next suppose $\mu _{i-1} = \de$. Then $2\de \le \pi$.
Since $\de \ne \pi/2$, there are other tiles with a vertex at $F_{n-1}$.
However, we have $2\de +\al >\de +2\al >\de +\be =\pi$, which is impossible.

Thus the only possibility is $\mu _{n-1} = \be$. Since  
$|F_{n-2} F_{n-1} |$ is $b$ or $d$, we have $|F_{n-2} F_{n-1}|=b$ and
$\la _{n-1} =\ga$. Then $\mu_{n-2} =\al$, $|F_{n-3} F_{n-2}|=d$
and $\mu _{n-2} =\de$. Continuing this argument we obtain 
$|F_{n-i-1} F_{n-i}|=d$, $\la _{n-i}=\de$, $\mu _{n-i}=\al$ if $i$ is even, and
$|F_{n-i-1} F_{n-i}|=b$, $\la _{n-i}=\ga$, $\mu _{n-i}=\be$ if $i$ is odd.

Thus $\la _1 = \ga$ or $\de$. Since the angle of $Q$ at $D$ is $\de$, we
have $\de +\la _1 \le \pi$. Then $\la _1 =\de$ and $2\de \le \pi$. Then,
by $\de \ne \pi/2$, there must be other tiles with a vertex at $D=F_0$.
However, as we saw above, we have $2\de +\al >\pi$, which is impossible.

This contradiction proves $\de =\be =\pi/2$. Then, by $\be <2\al$ we get
$\pi /4 <\al <\pi /2$.

The hypotenuse of the right triangles $ABC$ and $ACD$ is $AC$. Therefore, we have $a<\ol{AC}<c+d$, and thus 
\begin{equation}\label{e3a}
a-c<d.
\end{equation}
In the sequel we assume that none of (i)-(iii) of the lemma holds.
In particular, we have $\al \ne \pi /3$ and $\ga \ne 2\al$. Then $\ga$ is
not the linear combination of $\al$ and $\pi /2$ with nonnegative
integer coefficients. Indeed, this follows from $\ga =\pi -\al <3\pi /4 <
\al +\pi /2 <3\al$ and from $\ga \ne 2\al$.

We have $\la _1 =\al$, $\pi /2$ or $\ga$. Since the angle of $Q$ at $D$
equals $\de =\pi /2$, we have $\pi /2 +\la _1 \le \pi$, and thus $\la _1 =\al$
or $\pi /2$. Now $\la _1 =\al$ is impossible, as $\al +\pi /2 <\pi <2\al +\pi /2$. Thus $\la _1 =\pi /2$ and $\mu _1 =\al$ or $\ga$. If $\mu _1 =\ga$, then
$\la _2 =\al$. The converse is also true: if $\mu _{1} =\al$, then
$\la _2 =\ga$, since $\ga$ is not the linear combination of $\al$ and
$\pi /2$ with nonnegative integer coefficients. Therefore, we have
$\mu _1 =\ga$ and $\la _2 =\al$, or $\mu _1 =\al$ and $\la _2 =\ga$.

Let the vertices of $R_1$ and $R_2$ be $D=F_0 , F_1 ,G_1 ,H_1$, and
$F_1 , F_2 ,G_2 ,H_2$, respectively. Since $\la _1 =\al$ or $\ga$,
the angle of $R_1$ at $H_1$ equals $\ga$ or $\al$. Therefore, the point
$H_1$ is not a vertex of $Q'$. Indeed, $A$ is a vertex of $Q'$, and if $H_1$
was also a vertex, then the angle of $Q'$ at $H_1$ would be $\pi /2$, which
is impossible by $\ga >\pi /2$ and $\al <\pi /2 <2\al$. Thus there is a tile
$P$ having a vertex at $H_1$ and having angle $\mu _1$ at $H_1$. The point $H_1$
is the endpoint of two sides of $P$, one of them lies on the boundary of $Q'$,
the other one, $H_1 J$ lies in the interior of $Q'$ (apart from the point
$H_1$).

\begin{figure}
\includegraphics[width=3in]{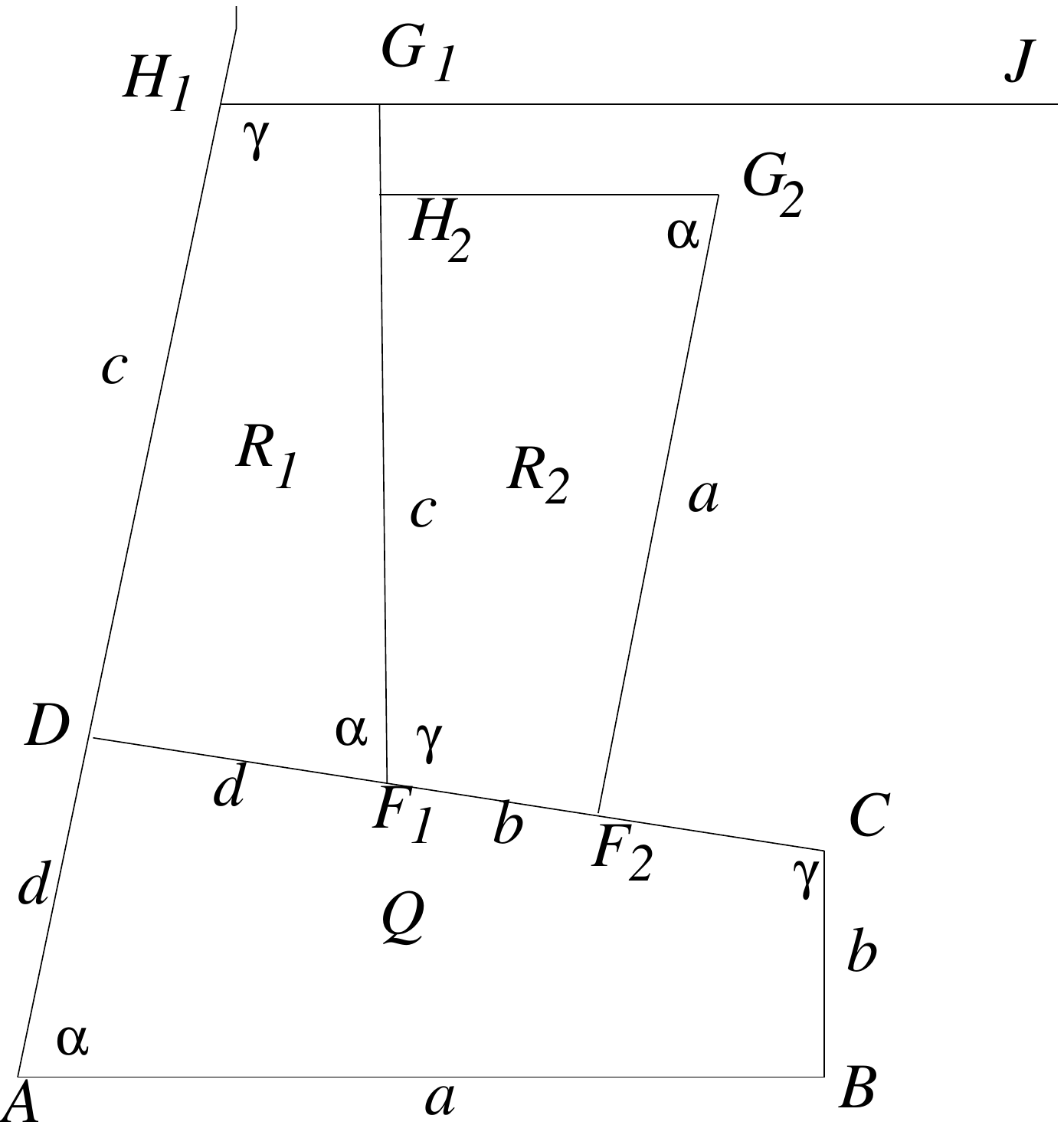}
\caption{}
\label{f2}
\end{figure}

\noi
{\bf Case III.2.1:} $\mu _1 =\al$ (see Figure \ref{f2}). Then $\ol{H_1 J}=a$
or $d$. Since $a> d>b$ by \eqref{e3}, the point $G_1$ is in the interior
of the segment $H_1 J$. The segments $G_1 J$ and $H_2 G_2$ are parallel to each
other and perpendicular to the segment $G_1 H_2$.

Then there is a partition $G_1 =K_0 \stb K_m  =H_2$
of the segment $G_1 H_2$ and there are tiles $P_1 \stb P_m$ such that 
$K_{i-1} K_i$ is a side of $P_i$ for every $i=1\stb m$. We have $m\ge 2$,
since the tile $P_1$ cannot have angle $\pi /2$ at two consecutive vertices. 

Now $\ol{G_1 H_2}=\ol{G_1 F_1} -\ol{H_2 F_1}=a-c$. By \eqref{e3a} we have
$a-c<d<c$, and thus $|K_{i-1} K_i |=b$ for every $i=1\stb m$.

Let $\rho _i$ and $\si _i$ denote the angle of $P_i$ at $K_{i-1}$ and at $K_i$,
respectively. Since $\rho _1 =\pi /2$, we have $\si _1 =\ga$, and thus
$\rho _2 =\al$. However, this implies $|K_{1} K_2 |\ne b$, which is impossible.

\begin{figure}
\includegraphics[width=3in]{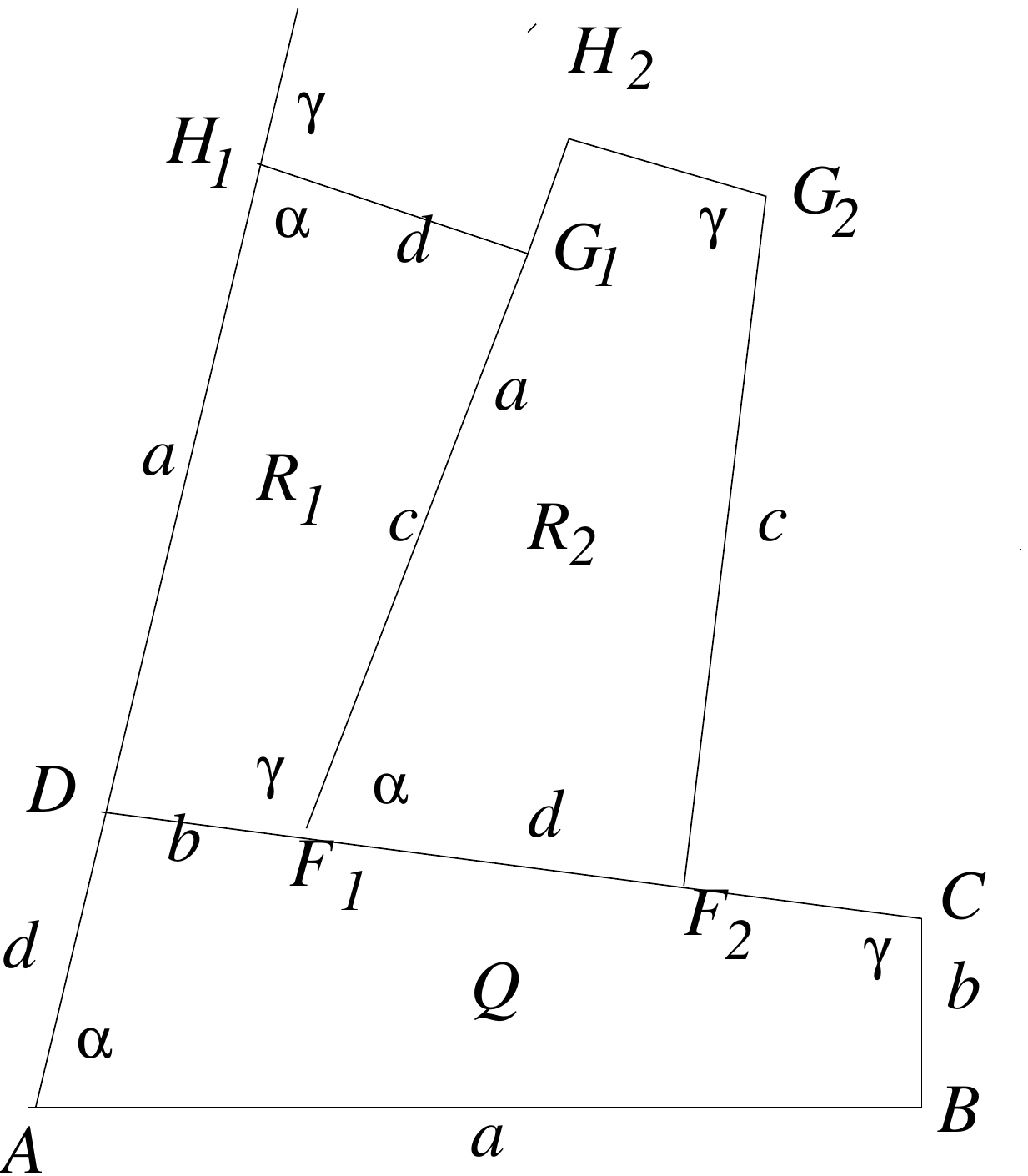}
\caption{}
\label{f3}
\end{figure}

\noi
{\bf Case III.2.2:} $\mu _1 =\ga$ (see Figure \ref{f3}). In this case
$\ol{F_1 H_2}=a>c=\ol{F_1 G_1}$.
Then there is a partition $H_1 =L_0 \stb L_t  =G_1$
of the segment $H_1 G_1$ and there are tiles $S_1 \stb S_t$ such that 
$L_{i-1} L_i$ is a side of $S_i$ for every $i=1\stb t$. If $t=1$, then
$H_1 G_1$ is a side of length $d$ of the tile $S_1$, and then the angles
of $S_1$ at $H_1$ and $G_1$ are $\al$ and $\pi /2$ in some order.
However, since the angle of $R_1$ at $H_1$ equals $\al$ and $\ga$ is not the
linear combination of $\al$ and $\pi /2$ with nonnegative integer coefficients,
the angle of $S_1$ at $H_1$ must be $\ga$, which is impossible.

Therefore, we have $t\ge 2$, and thus $|L_{i-1} L_i |=b$ for every $i$
by \eqref{e3}. Then the angles of $S_i$ at $L_{i-1}$ and $L_i$ are $\ga$
and $\pi /2$ in some order. Since $2\ga >\ga +\pi /2>\pi$, this is possible
only if $t=2$ and $S_2$ has angle $\ga$ at $G_1$. However, $S_2$ must have
a right angle at $G_1$, which is a contradiction. \hfill $\square$

\section{Ruling out families (i) and (ii)}
\begin{lemma}\label{l4}
Suppose $Q$ satisfies the conditions of {\rm (i)} of Theorem \ref{t2}. Then $Q$
is not a reptile.
\end{lemma}

\begin{figure}
\includegraphics[width=3in]{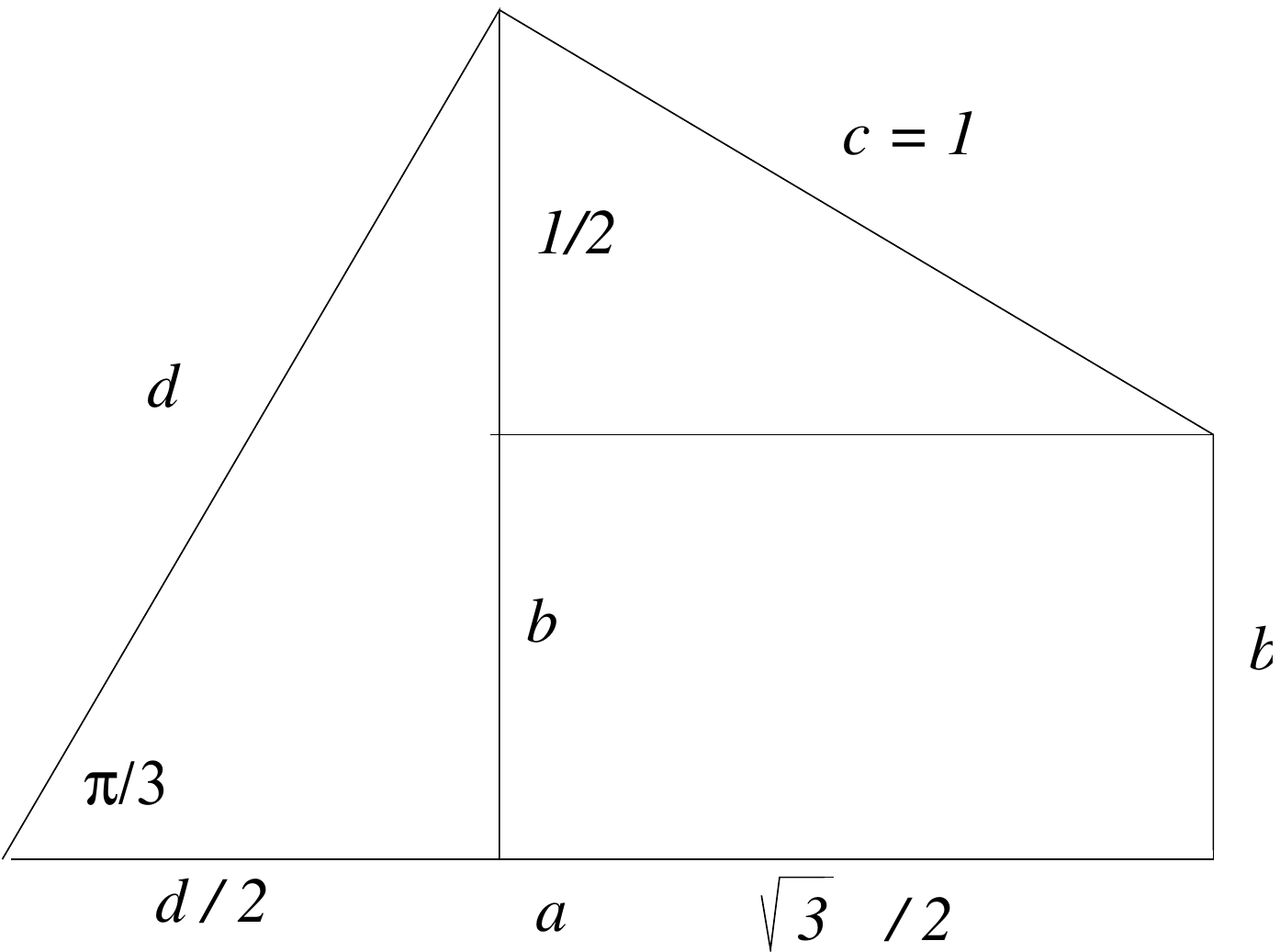}
\caption{}
\label{f4}
\end{figure}

\proof We may assume $c=1$. Suppose $Q$ is a reptile. Then there is a
quadrilateral $Q'$ similar to $Q$ such that $Q'$ can be tiled with $k\ge 2$
congruent copies of $Q$. Let $Q''$ be similar to $Q$ such that $Q''$ is tiled
with $k$ congruent copies of $Q'$. Then $Q''$ can be tiled with $k^2$ congruent
copies of $Q$, and thus the sides of $Q''$ have length $ka, kb , k$ and $kd$.
We have
\begin{equation}\label{e4}
  b=d\frac{\sqrt 3}{2} -\frac{1}{2} \ \text{and} \  a=\frac{d}{2} +
  \frac{\sqrt 3}{2} .
\end{equation}
(See Figure \ref{f4}.) Let $X,Y$ be consecutive vertices of $Q''$ such that
$\ol{XY} =kb$, and $Q''$ has a right angle at $X$. Let $p,q,r,s$ be the number
of tiles having a side of length $a,b,1,d$ lying on $XY$.

We show that $\max (p,s)>0$. Let $X=X_0 \stb X_n =Y$ be a partition of $XY$ such
that each $X_{i-1}X_i$ is a side of the tile $T_i$ $(i=1\stb n)$. If
$\ol{X_0 X_1}=a$ or $d$, then we have $\max (p,s)\ge 1$. If, however,
$\ol{X_0 X_1}=b$ or $1$ then the angle of $T_1$ at $X_1$ equals $\ga$. Then the
angle of $T_2$ at $X_1$ equals $\al$, and thus $\ol{X_0 X_1}=a$ or $d$, proving
$\max (p,s)>0$ in this case as well.

We have $kb=pa+qb+r+sd$, and thus $2(k-q)b=2pa+2r+2sd$. From $\max (p,s)>0$ it
follows that $k-q>0$. Substituting the values given by \eqref{e4} we get
$$(k-q)d\sqrt 3 -(k-q)=pd+p\sqrt 3 +2r +2sd.$$
Therefore, we have
$$[(k-q)\sqrt 3 -(p+2s)]d=p\sqrt 3 +(k-q+2r).$$
Multiplying both sides by $(k-q)\sqrt 3 +(p+2s)$ we obtain $Ad=B+C\sqrt 3$,
where $A,B,C\in \zz$, and $B>0$, $C>0$.
Indeed, $B=3(k-q)p +(p+2s)(k-q+2r)>0$ and $C=(k-q)(k-q+2r)+(p+2s)p>0$
by $k-q>0$ and $\max (p,s)>0$. Since $d>0$, it follows that $A>0$.
We find that $d=x+y\sqrt 3$, where $x,y$ are positive rational numbers.

Let $ZY$ be the side of $Q''$ of length $k$, and let $t,u,v,w$ be the number of
tiles having a side of length $a,b,1,d$ lying on $ZY$. One can show, similarly
to the previous  case, that $\max (t,w)>0$. We have $k-v=ta+ub+wd$, and thus
\begin{align*}
2(k-v)&=2ta+2ub+2wd=\\
&=  (tx+ ty\sqrt 3 + t\sqrt 3 ) +(ux\sqrt 3 +3uy-u)+(2wx +2wy\sqrt 3 ).
\end{align*}
Since $\sqrt 3$ is irrational, we have $t=w=0$, which is
impossible. \hfill $\square$

\begin{lemma}\label{l5}
Suppose $Q$ satisfies the conditions of {\rm (ii)} of Theorem \ref{t2}.
If $Q$ is not a trapezoid, then $Q$ is not a reptile.
\end{lemma}

\proof We may assume that $b=c=1$. Then
the triangle $BCD$ is isosceles, the length of its legs is $1$ and the angle
at the apex is $2\pi /3$. Therefore, the circumscribed circle $K$ of the
triangle $BCD$ has radius $1$. Since $K$ is the circumscribed circle of $Q$ as
well, we have $2\ge a\ge d>1$. Now $a=2$ happens only when $a$ is the diameter
of $K$, and then $Q$ is a trapezoid. Therefore, we have $2>a\ge d>1$,
It is easy to check that if $a=d$, then $\be =\de =\pi /2$ when, by
Lemma \ref{l4}, $Q$ is not a reptile. So we may assume that $2>a>d>1$,

Since $Q$ is not a trapezoid, $\al <\be <\ga$ and $\al <\de <\ga$.

Suppose $Q$ is a reptile. Then there is a quadrilateral $Q'$ similar to $Q$
such that $Q'$ can be tiled with $k\ge 2$ congruent copies of $Q$. 
The angle $\be$ of $Q'$ is packed with one single tile, since $\al <\be$
and $2\al =2\pi /3 =\ga >\be$. We may assume that
$Q$ itself is this tile; that is, $B$ is a vertex of $Q'$, and the sides
$AB$ and $BC$ lie on consecutive sides of $Q'$.

There is a tile $P$ having a vertex at $C$ and having angle $\al$ at $C$.
The point $C$ is the endpoint of two sides of $P$, one of them lies on the
boundary of $Q'$, the other one, $CE$ lies in the interior of $Q'$ (apart
from the point $C$). The length of $CE$ is $a$ or $d$, and both are greater
than $\ol{BC}=1$ by Lemma \ref{l1}.

Then it follows that there is a partition $A =A_0 \stb A_n  =D$
of the segment $AD$ and there are tiles $P_1 \stb P_n$ such that 
$A_{i-1} A_i$ is a side of $P_i$ for every $i=1\stb n$.
Since $\ol{AD}=d$ and $2> a> d>1$, we must have $n=1$; that is, $Q$ shares its
side $AD$ with another tile $R$.

Let $\la$ be the angle of $R$ at the vertex $D$. Since $\ol{AD}=d$, we have
$\la =\al$ or $\de$. The angle of $Q$ at $D$ is $\de$. If $\la =\al$, then
$\de +\al <\de +\be =\pi$ implies that there are other tiles with a vertex at
$D$, and thus $\de +2\al \le \pi$. However, $\de >\al =\pi -2\al$, so this case
is impossible.

If $\la =\de$, then $2\de \le \pi$. If $2\de =\pi$, then $\de =\be =\pi /2$.
In this case $Q$ is not a reptile by Lemma \ref{l4}. If $2\de <\pi$, then
there are other tiles with a vertex at $D$, and thus $2\de +\al \le \pi$.
However, $2\de +\al >\de +2\al >\pi$, so this case
is also impossible. \hfill $\square$

\section{Ruling out family $\ff _3$: preliminaries} In the next two sections our
aim is to prove the following.
\begin{lemma}\label{l6}
Suppose $Q$ satisfies the conditions of {\rm (iii)} of Theorem \ref{t2}.
If $Q$ is not a square, then $Q$ is not a reptile.
\end{lemma}
We may assume that $c=d=1$. Then we have $a>1>b$. From $a^2 +b^2 =2$ we get
$a<\sqrt 2$. Also, we have $a+b<2$, since $(a+b)/2<\sqrt{(a^2 +b^2 )/2}=1$.

If $\al =\pi /3$, then $a=1+b$ by \eqref{e4}. The converse is also true,
since the system of equations $a=1+b$, $a^2 +b^2 =2$ determines $a$ and $b$.
If $\al =\pi /3$, then $Q$ is not a reptile by Lemma \ref{l4}. Therefore, in
the sequel we may assume that $\al \ne \pi /3$ and $a\ne 1+b$.

Note that $\be = \pi /2 <2\al$ by Lemma \ref{l2}, and thus $\pi /4 <\al <
\pi /2$. Then $\ga$ is not the linear combination of $\al$ and $\pi /2$ with
nonnegative integer coefficients. Indeed, this follows from the inequalities
$3\al > \pi -\al = \ga$, $\al + \pi /2 > 3\pi /4 >\ga$ and from
$2\al \ne \pi -\al =\ga$.

Next we rule out two other special cases.

Suppose that $b=1/2$. From $a^2 +b^2 =2$ we get $a=\sqrt 7 /2$.
Suppose $Q$ is a reptile. Then there is a quadrilateral $Q''$ similar to $Q$
such that $Q''$ can be tiled with $k^2 \ge 4$ congruent copies of $Q$.
Let the vertices of $Q''$ be $X,Y,V,W$ such that $\ol{XY}=ka$, $\ol{YV}=kb$,
$\ol{VW}=\ol{WX}=k$. The side $XY$ of $Q''$ is packed with sides of some tiles.
Therefore, we have $ka=pa+qb+r$, where $p,q,r$ are nonnegative integers.
Since $b$ is rational and $a$ is irrational, we have $q=r=0$ and $k=p$; that is,
the side $XY$ is only packed with sides of tiles of length $a$. Let
$X=X_0 , X_1 \stb X_k =Y$ be a partition of the segment $XY$ and let
$T_i$ $(i=1\stb k)$ be tiles such that $X_{i-1} X_i$ is a side of length $a$ of
the tile $T_i$ for every $i=1\stb k$. Since the angle of $T_k$ at $Y$
equals $\pi /2$, it follows that the angle of $T_k$ at $X_{k-1}$ is $\al$.
Then the angle of $T_{k-1}$ at $X_{k-1}$ equals $\ga$, but then $\ol{X_{k-2} X_k}$
equals either $b$ or $1$, which is impossible.

Next suppose that $a=2-2b$. Then $a^2 +b^2 =2$ gives $a=(2 +2\sqrt 6 )/5$ and
$b=(4-\sqrt 6 )/5$.

Suppose $Q$ is a reptile. Then there is a quadrilateral $Q''$ similar to $Q$
such that $Q''$ can be tiled with $k^2 \ge 4$ congruent copies of $Q$.
Let the vertices of $Q''$ be $X,Y,V,W$ such that $\ol{XY}=ka$, $\ol{YV}=kb$,
$\ol{VW}=\ol{WX}=k$. The side $YV$ of $Q''$ is packed with sides of some tiles.
Therefore, we have $kb=pa+qb+r$, where $p,q,r$ are nonnegative integers.
Then $5(k-q)b=5pa+5r$; that is,
$$(k-q)(4-\sqrt 6 )=p(2+2\sqrt 6 )+5r.$$
Since $\sqrt 6$ is irrational and $k-q, p,r$ are nonnegative integers,
we have $k-q=p=r=0$. That is, the side $YV$ is only packed with
sides of tiles of length $b$. However, each of these tiles has angle
$\ga$ at one of the endpoints of its side of length $b$, which makes such a
tiling impossible.

Therefore, we also assume that $b\ne 1/2$ and $a\ne 2-2b$.

Suppose a convex subset $\Si \su \sik$ is tiled with congruent copies of $Q$.
We denote by $U$ the union of the boundaries of the tiles. By a {\it barrier}
we mean
a broken line $X_1 X_2 X_3 X_4$ covered by $U$ such that the segments $X_1 X_2$
and $X_3 X_4$ lie on the same side of the line going through $X_2 X_3$. We say
that $X_2 X_3$ is the {\it base} of the barrier, and the angles
$X_1 X_2 X_3 \szog$ and $X_2 X_3 X_4 \szog$ are the {\it angles} of the barrier.
By the {\it halfplane} of the barrier we mean the halfplane bounded by the
line going through $X_2 X_3$ and containing the segments $X_1 X_2$ and $X_3 X_4$.

If $X_1 X_2 X_3 X_4$ is a barrier, then there are tiles $T_1 \stb T_n$ and
there is a partition $X_2 =Y_0 ,Y_1 \stb Y_n =X_3$ of the base of the barrier
such that $Y_{i-1} Y_i$ is a side of $T_i$ $(i=1\stb n)$. We say that
$T_1 \stb T_n$ is a {\it tiling} of the barrier, and 
$X_2 =Y_0 ,Y_1 \stb Y_n =X_3$ is the partition of the barrier corresponding
to the tiling.

\begin{figure}
\includegraphics[width=3in]{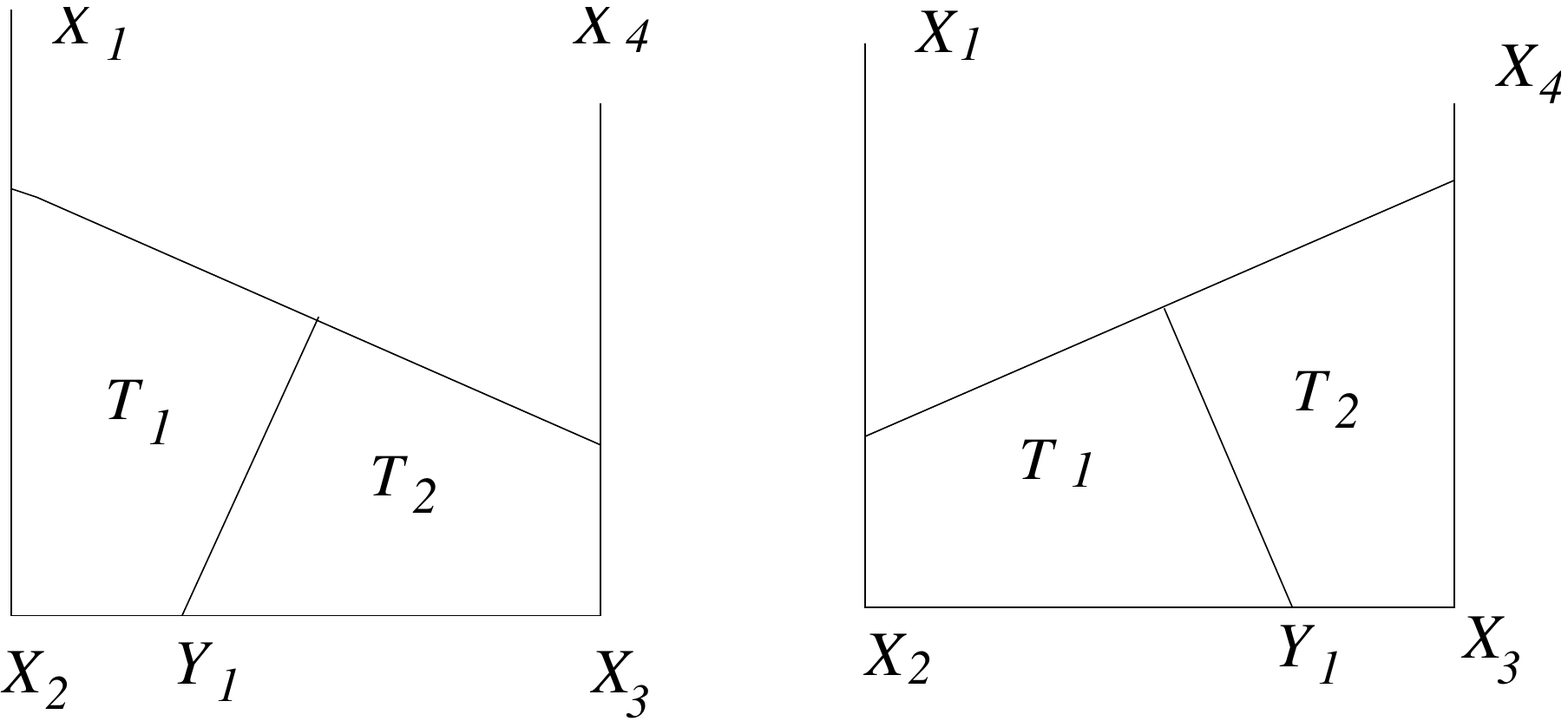}
\caption{}
\label{f8b}
\end{figure}

\begin{lemma}\label{l8}
Let $X_1 X_2 X_3 X_4$ be a barrier with right angles.
\begin{enumerate}[{\rm (i)}]
\item If $T_1 \stb T_n$ is a tiling of the barrier and $X_2 =Y_0 ,Y_1 \stb
Y_n =X_3$ is the partition corresponding to the tiling, then $n$ is even, and
$\ol{Y_{2i -2}Y_{2i}}$ equals one of $1+b, a+b, 2, 1+a$ for every $i=1\stb n/2$. 
\item We have $\ol{X_2 X_3}>1$.
\item If $\ol{X_2 X_3}=1+b$, then $\min(\ol{X_1 X_2} ,\ol{X_3 X_4} )\le 1$. 
\item If $\ol{X_2 X_3}=a+b$, then there are only two possible tilings as
shown by Figure \ref{f8b}.
\end{enumerate}
\end{lemma}

\proof (i) Since $T_1$ has a right angle at $Y_0$, its angle at $Y_1$ is
$\al$ or $\ga$. If it is $\ga$, then the angle of $T_2$ at $Y_1$ is $\al$. Then
$\ol{Y_{0}Y_{1}} \in \{ 1,b\}$, $\ol{Y_{1}Y_{2}} \in \{ 1,a\}$ and
$\ol{Y_{0}Y_{2}} \in \{ 1+b ,a+b, 2, 1+a\}$.

If the angle at $Y_1$ is $\al$, then the angle of $T_2$ at $Y_1$ is $\ga$,
because $\ga$ is not the linear combination of $\al$ and $\pi /2$ with
nonnegative integer coefficients. Then we obtain $\ol{Y_{0}Y_{2}} \in
\{ 1+b ,a+b, 2,1+a\}$ again. Since the angle of $T_2$ at $Y_2$
is $\pi /2$, $Y_2 Y_n$ is the base of a barrier with right angles, and we
can argue by induction. 

\noi
(ii): This follows from (i).

\begin{figure}
\includegraphics[width=2.6in]{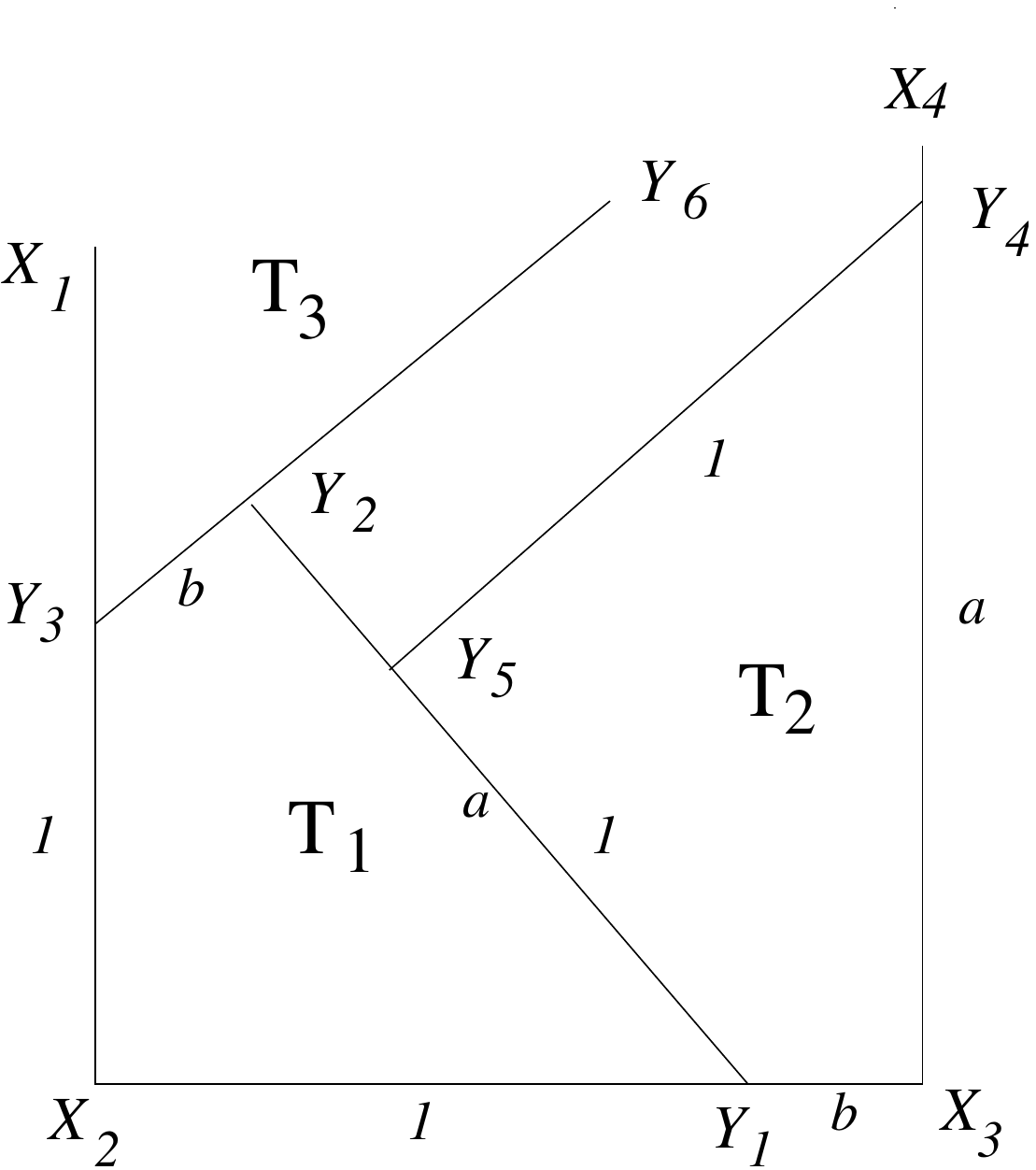}
\caption{}
\label{f8a}
\end{figure}

\noi
(iii) Suppose $\min(\ol{X_1 X_2} ,\ol{X_3 X_4} )> 1$. There are two tilings
of the barrier. One is shown by Figure \ref{f8a}, the other is obtained by
reflecting this tiling about the perpendicular bisector of the base. By
symmetry, we may assume that the tiling is as in Figure \ref{f8a}.

Since $\ol{X_1 X_2}>1$, there is a tile $T_3$ having the point $Y_3$ as vertex,
and having angle $\al$ at $Y_3$. The tile $T_3$ has a side $Y_3 Y_6$ of
length $1$ or $a$, and thus the point $Y_2$ is an inner point of the segment
$Y_3 Y_6$. Then $Y_6 Y_2 Y_5 Y_4$ is a barrier with two right angles
and base length $a-1<1$, which is impossible by (ii).

\noi
(iv) Let $T_1 \stb T_n$ be a tiling, and let $X_2 =Y_0 ,Y_1 \stb Y_n =X_3$ be
the corresponding partition. Since $a+b <2$, we have either $\ol{Y_0 Y_1}=b$
or $\ol{Y_{n-1} Y_n}=b$. Suppose the former. Then $T_1$ has angle $\ga$ at
$Y_1$, and thus $T_2$ has angle $\al$ at $Y_1$. If $\ol{Y_1 Y_2}=1$, then
$Y_3 \stb Y_n$ constitute a tiling of a barrier with two right angles and base
length $a-1<1$, contradicting (ii). Therefore, we have $\ol{Y_1 Y_2}=a$. Then
$n=2$, and the tiling looks like the first tiling in Figure \ref{f8b}.
We obtain the second, if $\ol{Y_{n-1} Y_n}=b$. \hfill $\square$

\begin{lemma}\label{l9}
Let $X_1 X_2 X_3 X_4$ be a barrier such that $X_1 X_2 X_3 \szog =\pi /2$,
$X_2 X_3 X_4 \szog =\ga$, and $\ol{X_2 X_3}=1$. Then the tiling of the barrier
consists of one single tile having $X_2 X_3$ as a side.
\end{lemma}

\proof Let $T_1 \stb T_n$ be a tiling of the barrier, and let $X_2 =Y_0 ,Y_1
\stb Y_n =X_3$ be the partition corresponding to the tiling. Since the angle
of $T_n$ at $X_3$ equals $\ga$, the angle of $T_n$ at $Y_{n-1}$ equals $\pi /2$.
If $Y_{n-1} \ne X_2$, then $\ol{Y_{n-1} Y_n}=b$, and $1-b=\ol{X_2 Y_{n-1}}$
is the base length of a barrier with two right angles. This, however,
contradicts (ii) of Lemma \ref{l8}. \hfill $\square$

\begin{lemma}\label{l10}
Let $X_1 X_2 X_3 X_4$ be a barrier such that $X_1 X_2 X_3 \szog =\pi /2$,
$X_2 X_3 X_4 \szog =\ga$, and $\ol{X_2 X_3}=a$. Then $a=1+2b$, 
the tiling of the barrier is unique, and consists of three tiles $T_1 ,T_2 ,T_3$
such that $\ol{Y_0 Y_1}=1$ and $\ol{Y_1 Y_2}=\ol{Y_2 Y_3}=b$, where
$X_2 =Y_0 ,Y_1 , Y_2 ,  Y_3 =X_3$ is the partition corresponding to the tiling.
\end{lemma}

\proof Let $T_1 \stb T_n$ be a tiling of the barrier, and let
$X_2 =Y_0 ,Y_1 \stb Y_n =X_3$ be the partition corresponding
to the tiling. Since the angle of $T_n$ at $X_3$ equals $\ga$, the angle
of $T_n$ at $Y_{n-1}$ equals $\pi /2$. Therefore, $X_2 Y_{n-1}$ is the base
of a barrier with two right angles. By Lemma \ref{l8}, $\ol{X_2 Y_{n-1}}
>1$. Then $\ol{Y_{n-1} X_3}= \ol{X_2 X_3}- \ol{X_2 Y_{n-1}}<a-1<1$, and 
thus $\ol{Y_{n-1} X_3}=b$ and $\ol{X_2 Y_{n-1}}=a-b$. Since $a-b<2(1+b)$, we have
$a-b\in \{ 1+b ,a+b ,2, 1+a\}$ by (i) of Lemma \ref{l8}. We obtain
$a-b =1+b$, $a=1+2b$ and $n=3$.

\begin{figure}
\includegraphics[width=3in]{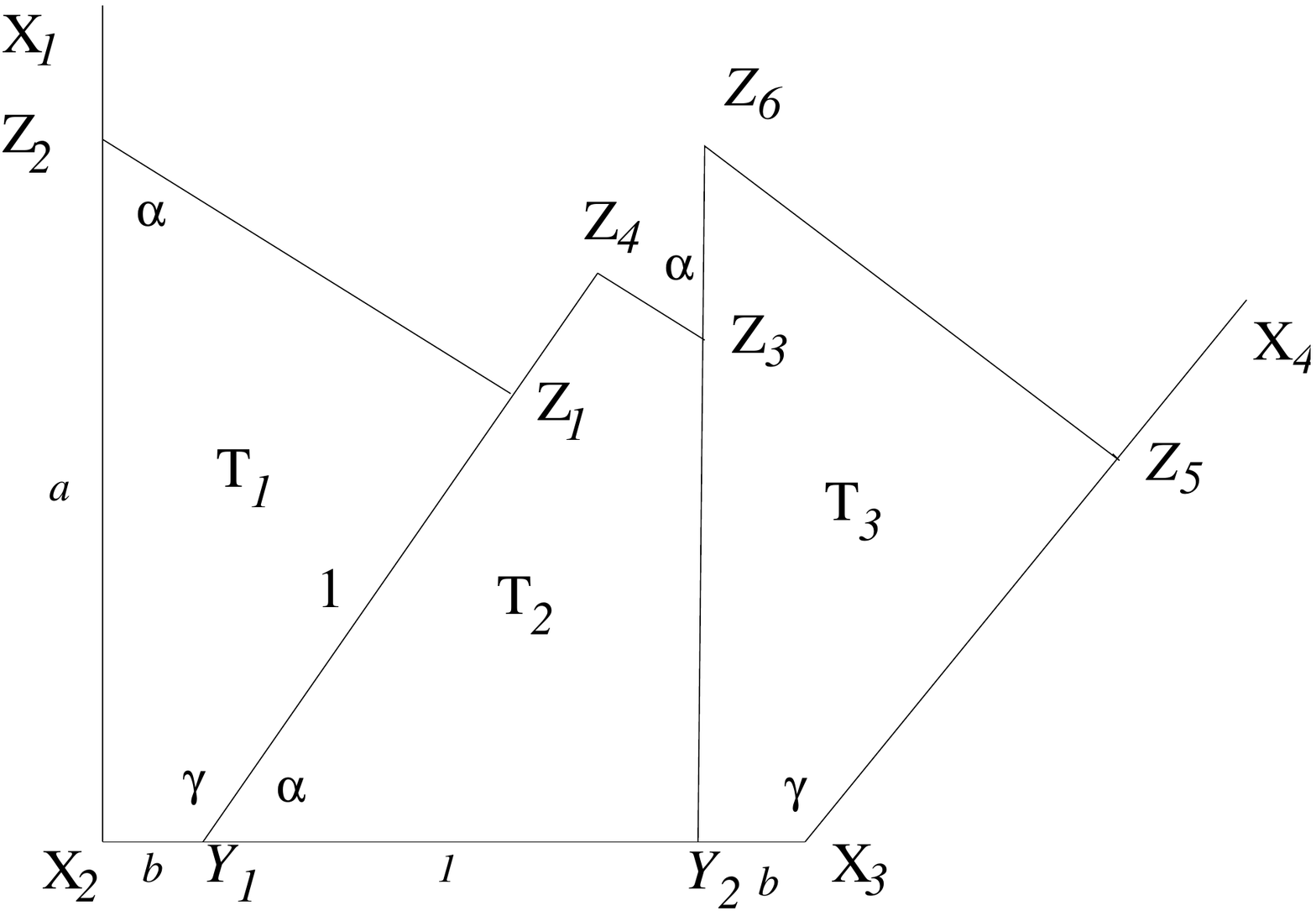}
\caption{}
\label{f6}
\end{figure}

In order to complete the proof we have to show that $\ol{Y_0 Y_1}=1$ and
$\ol{Y_1 Y_2}=b$. If this is not true, then $\ol{Y_0 Y_1}=b$ and
$\ol{Y_1 Y_2}=1$. Figure \ref{f6} shows the arrangement of the tiles in this
case. Then $Z_2 X_2 Y_2 Z_6$ is a barrier with two right angles and base length
$1+b$. Since $\ol{Z_2 X_2} =a>1$ and $\ol{Z_6 Y_2}=a>1$, this contradicts
(iii) of Lemma \ref{l8}. \hfill $\square$

\begin{figure}
\includegraphics[width=3in]{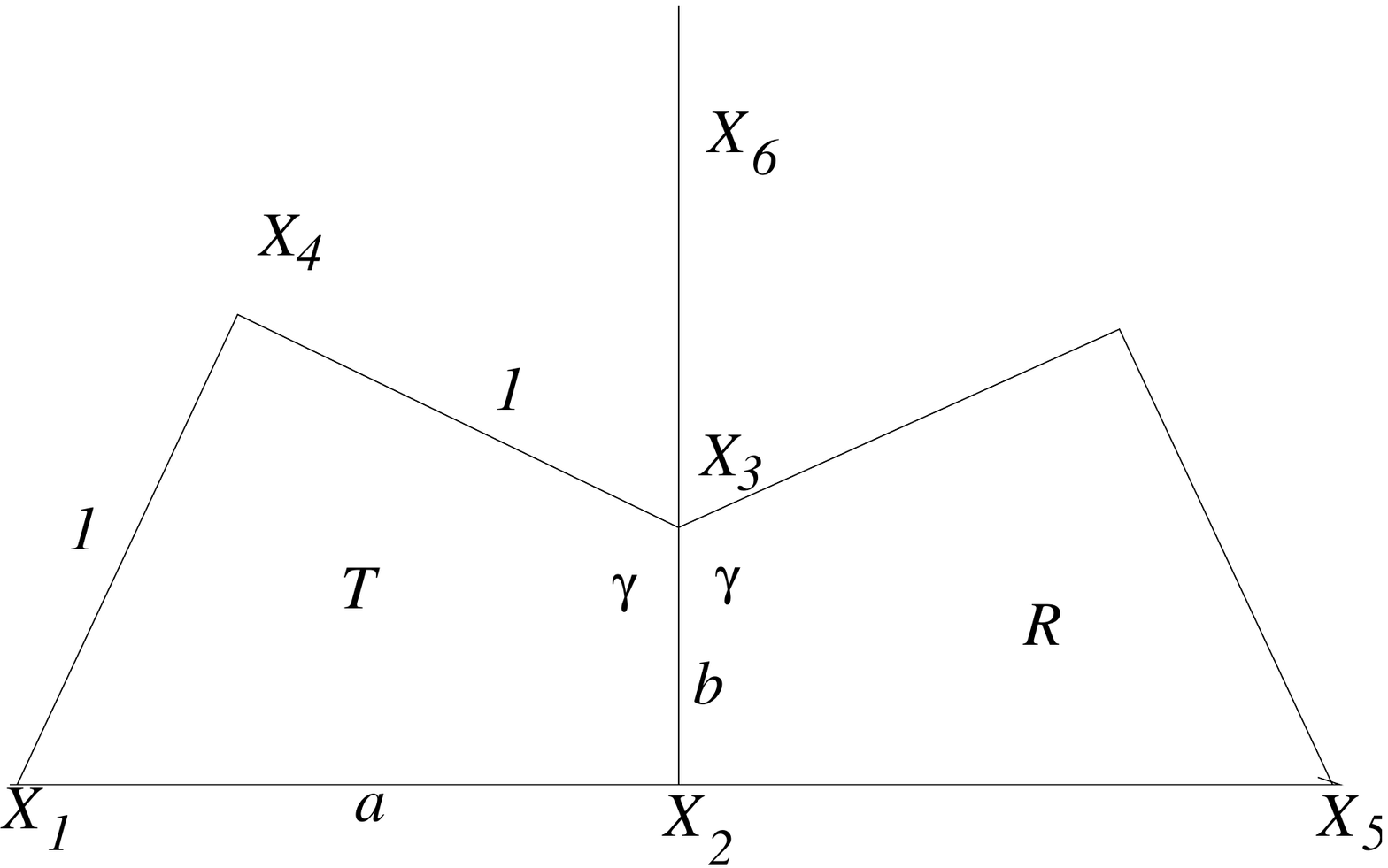}
\caption{}
\label{f5}
\end{figure}

\begin{lemma}\label{l11}
Let $T$ be a tile with vertices $X_1 ,X_2 , X_3 ,X_4$ such that $\ol {X_2 X_3}=b$
and $\ol {X_3 X_4}=1$ (see Figure \ref{f5}). Suppose that there is a segment
$X_1 X_5 \su U$ such that $X_2$ is an inner point of $X_1 X_5$. Then there is a
segment $X_2 X_6 \su U$ such that $X_3$ is an inner point of $X_2 X_6$.
\end{lemma}

\proof Let $R$ be the tile having a vertex at $X_2$, different from $T$
and lying on the same side of $X_1 X_5$ as $T$. Then $R$ has a side $X_2 X_7$
lying on the line $\ell$  going through $X_2$ and $X_3$.

If $\ol{X_2 X_7}>b$, then we can take $X_6 =X_7$.
If $\ol{X_2 X_7}=b$, then $X_7 =X_3$, and the angle of $R$ at $X_3$ is $\ga$.
Thus $T$ and $R$ both have an angle $\ga$ at $X_3$. Since $2\pi -2\ga
=2\al$, it follows that there are two other tiles having a vertex at $X_3$,
and both tiles have an angle $\al$ at $X_4$. Since $\ga +\al =\pi$, we can see
that the segment $X_2 X_3$ can be continued in $U$; that is, there is a point
$X_6$ on $\ell$ such that $X_3$ is an inner point of a segment
$X_2 X_6 \su U$. \hfill $\square$

\begin{lemma}\label{l12}
Let the segments $XY$ and $YZ$ be be covered by $U$, and suppose $XYZ\szog =
\pi /2$ and $\ol {XY}>1$, $\ol {YZ}>1$. Let $\ell _X$ and $\ell _Z$ denote the
halflines that start from $Y$ and go through $X$ and $Z$, respectively.
Let $T$ be the (unique) tile having a vertex at $Y$ and lying in the quadrant
bounded by the halflines $\ell _X$ and $\ell _Y$. Then the sides of $T$
with endpoint $Y$ have lengths $a$ and $b$.
\end{lemma}

\begin{figure}
\includegraphics[width=4in]{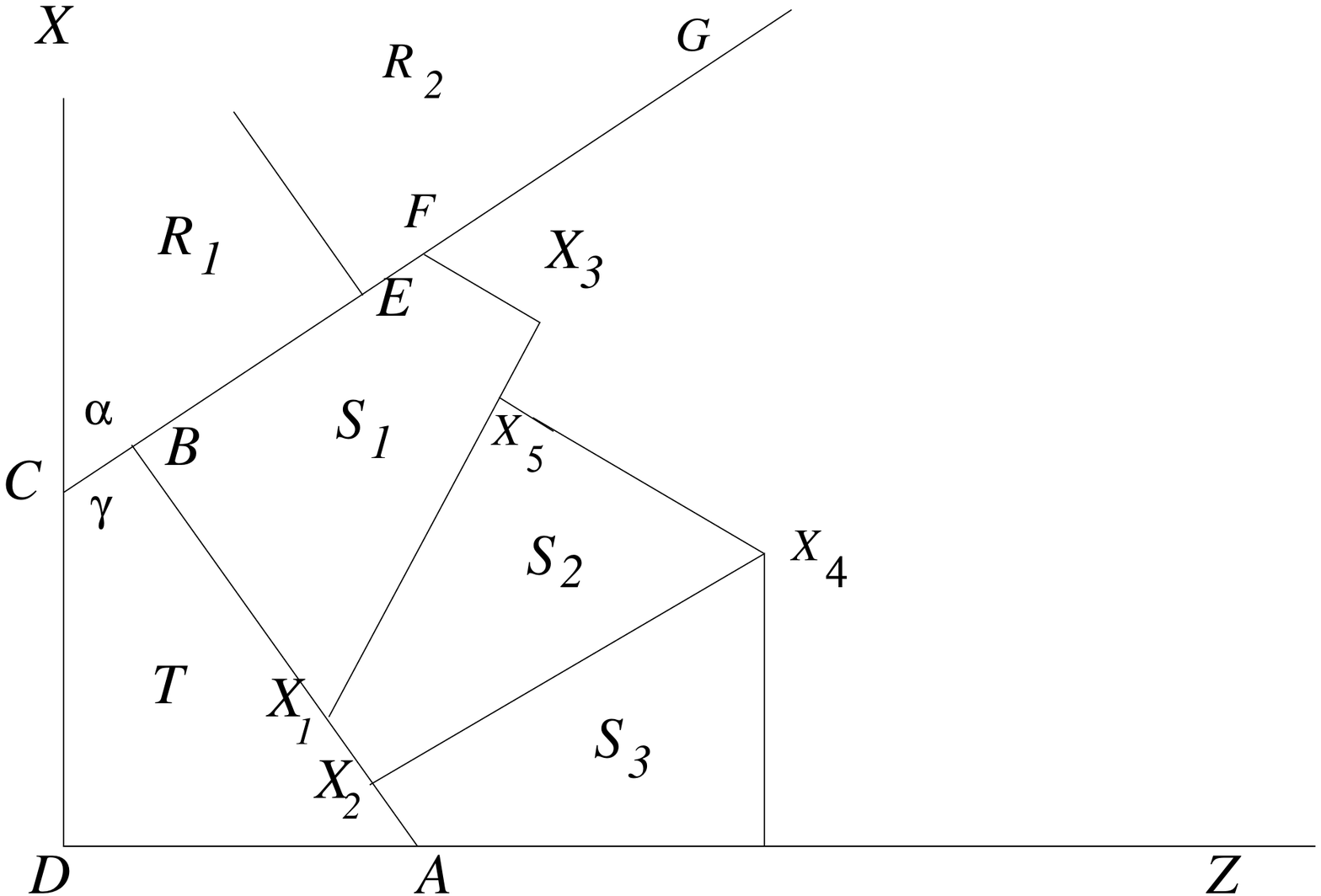}
\caption{}
\label{f7}
\end{figure}

\proof Since $XY, YZ\su U$ and $\al >\pi /4$, it follows that $T$ has a
right angle at $Y$. If the statement of the lemma is not true, then both
sides of $T$ with endpoint $Y$ have length $1$. We may assume that $T=Q$
with the usual labeling. Then $Y=D$, and by symmetry we may also assume that
$A$ is an inner point of the segment $DZ$ and $C$ is an inner point of the
segment $DX$ (see Figure \ref{f7}).

Since $\ol{XD}>1$, there is a tile $R_1$ having a vertex at $C$ and having
angle $\al$ at $C$, and such that the point $B$ is an inner point of the side
$CE$ of $R_1$. Note that $\ol{CE}$ equals $1$ or $a$, and that $R_1$ has a right
angle at $E$.

Now $EBAZ$ is a barrier with angle $\pi /2$ and $\ga$ and with base length $a$.
By Lemma \ref{l10}, we have $a=1+2b$, and the tiling of $EBAZ$ is unique, and
consists of three tiles $S_1 ,S_2 ,S_3$ as in Figure \ref{f7}. Let $B,X_1 ,X_2 ,
A$ be the partition corresponding to the tiling of $EBAZ$; then $\ol{BX_1}=1$
and $\ol{X_1 X_2}=\ol{X_2 A}=b$. Let $B,X_1 ,X_3 ,F$ be the vertices of $S_1$,
and let $X_1 , X_2 ,X_4 ,X_5$ be the vertices of $S_2$.

Let $\ell$ denote the line going through the points $C,B$ and $E$. We prove
that there is a point $G$ on $\ell$ such that $CG \su U$ and $F$ is an inner
point of the segment $CG$.

If $\ol{CE}=a=1+2b$, then $\ol{CE}>1+b =\ol{CF}$, and we can take $G=E$.

Suppose $\ol{CE}=1$, and let $R_2$ be the tile having a vertex at $E$ and
different from $R_1$. Then $R_2$ has a side $EH$ lying on $\ell$. If
$\ol{EH}>b$, then we can take $G=H$. If $\ol{EH}=b$, then $H=F$, and the
angle of $R_2$ at $F$ is $\ga$. Thus $R_2$ and $S_1$ both have an angle $\ga$ at
$F$. Since $2\pi -2\ga =2\al$, it follows that there are two other tiles having
a vertex at $F$, and both tiles have an angle $\al$ at $F$. Since
$\ga +\al =\pi$, we can see that the segment $CF$ can be continued in $U$; that
is, there is a point $G$ with the required properties. 

Then $GBX_2 X_4$ is a barrier having right angles and base length
$1+b$. Since $\ol{GB}>1$ and $\ol{X_2 X_4}=a>1$, this contradicts (iii) of
Lemma \ref {l8}. This contradiction completes the proof. \hfill $\square$

\begin{lemma}\label{l9a}
Let $X_1 X_2 X_3 X_4$ be a barrier such that $X_1 X_2 X_3 \szog =\pi /2$,
$X_2 X_3 X_4 \szog =\ga$, and $\ol{X_2 X_3}=2a+b$. Then $\ol{X_1 X_2} \le 1$.
\end{lemma}

\proof Let $T_1 \stb T_n$ be a tiling of the barrier, and let
$X_2 =Y_0 ,Y_1 \stb Y_n =X_3$ be the corresponding partition.

Since $T_n$ has angle $\ga$ at $X_3$, we have either $\ol{Y_{n-1} Y_n}=1$ or
$\ol{Y_{n-1} Y_n}=b$. We consider the two cases separately.

\noi
{\bf Case I:} $\ol{Y_{n-1} Y_n}=1$. The tiles $T_1 \stb T_{n-1}$ constitute a
tiling with two right angles and with
base length $2a+b-1$. By (i) of Lemma \ref{l8}, $n-1$ is even, and $2a+b-1$
equals a sum $c_1 +\ldots +c_k$, where $k=(n-1)/2$, and each $c_i$ equals
one of $1+b, a+b,2,1+a$. Now
$a<\sqrt 2$ implies $2a<3$ and $2a+b-1<2+b<2(1+b)$. Therefore, we have $k=1$
and $2a+b-1\in \{ 1+b, a+b,2,1+a \}$. By $a+b <2a+b-1<1+a$, the only
possibility is $2a+b-1=2$, and thus we have $\ol{Y_0 Y_1}=\ol{Y_1 Y_2}=
\ol{Y_2 Y_3}=1$.

Let $X_2 , Y_1 ,Z_1 , Z_2$ be the vertices of $T_1$. Then
$\ol{X_2 Y_1}=\ol{X_2 Z_1}=1$. Now $\ol{X_2 X_3} =2a+b>1$. Therefore,
$\ol{X_1 X_2}\le 1$ follows from Lemma \ref{l12}.

\noi
{\bf Case II:} $\ol{Y_{n-1} Y_n}=b$. If the vertices of $T_n$ are $Y_{n-1} , X_3 ,
W_1 ,W_2$, then $\ol{W_2 Y_{n-1}}=a$.
The tiles $T_1 \stb T_{n-1}$ constitute a tiling with two right angles and with
base length $2a$. By (i) of Lemma \ref{l8}, $n-1$ is even, and $2a$
equals a sum equals a sum $c_1 +\ldots +c_k$, where $k=(n-1)/2$, and each $c_i$
equals one of $1+b, a+b,2,1+a$. Now 
$2a>1+a$ implies that $k \ge 2$. On the other hand,
$a<\sqrt 2$ implies $2a<3$, hence $k \le 2$. Therefore, we have
$k = 2$, and $2a$ is the sum of two elements of the set $\{ 1+b, a+b,2,
1+a \}$. Since the case $a=1+b$ has been excluded and $2a<2(a+b)$, the only
possibility is $2a= (1+b)+(a+b)$. Then $\ol{Y_0 Y_2}=1+b$ and
$\ol{Y_2 Y_4}=a+b$, or the other way around.

\begin{figure}
\includegraphics[width=3.5in]{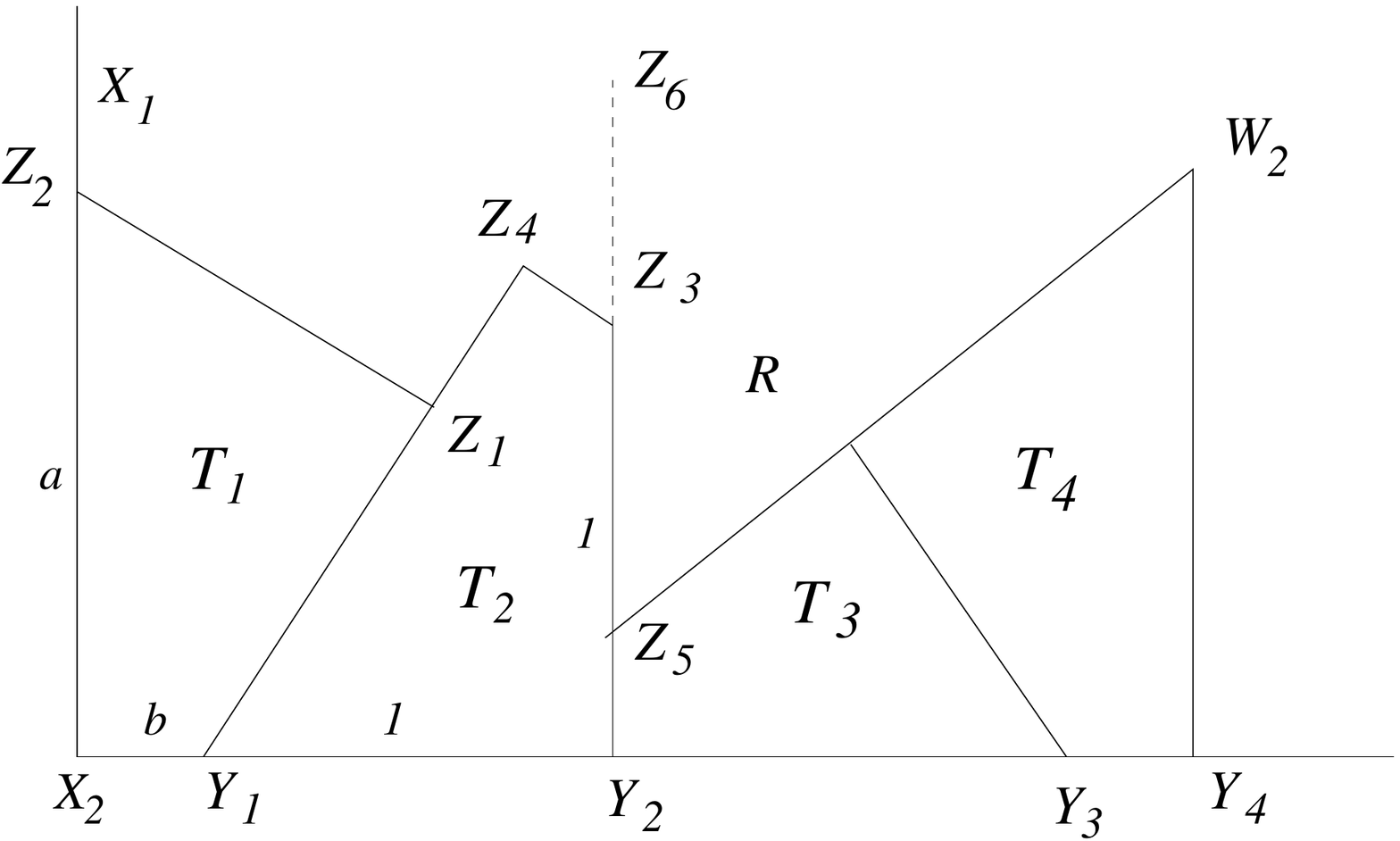}
\caption{}
\label{f8c}
\end{figure}

\noi
{\bf Case II.1:} $\ol{Y_0 Y_2}=1+b$ {\it and} $\ol{Y_2 Y_3}=a+b$.
Suppose $\ol{X_1 X_2}>1$. It follows from Lemma \ref{l12} that $\ol{Y_0 Y_1}=b$,
$\ol{Y_1 Y_2}=1$, and the arrangement of $T_1$ and $T_2$ is as in Figure
\ref{f8c}. By (iii) of Lemma \ref{l8},
the segment $Y_2 Z_3$ cannot be continued in $U$. Then, by (iv) of Lemma
\ref{l8}, the tiles $T_3$ and $T_4$ are as in Figure \ref{f8c}. Then there
is a tile $R$ having a vertex at $Z_5$. Then $R$ has a side $Z_5 Z_6$ such that
$Z_3$ is an inner point of $Z_5 Z_6$, contradicting (iii) of Lemma \ref{l8}.

\noi
{\bf Case II.2:} $\ol{Y_0 Y_2}=a+b$ {\it and}  $\ol{Y_2 Y_3}=1+b$.
The proof in this case is the same as in Case II.1. The only difference is
that we work from right to left; that is, from $Y_4$ to $X_2$,
using $\ol{W_2 Y_4 }=a>1$.  \hfill $\square$

\section{Ruling out family (iii): conclusion}
As we mentioned earlier, the square of side length $a+b$ can be tiled with four
copies of $Q$ (see Figure \ref{f0}).

Let $\Si$ denote the quadrant $\{ (x,y)\colon x\ge 0, \
y\ge 0\}$. Since $\Si$ is tiled with the squares
$$\si _{i,j} =[(i(a+b), (i+1)(a+b)]\times [j(a+b), (j+1)(a+b)]$$
$(i,j=0,1,\ldots )$, it follows that $\Si$ can be tiled with congruent
copies of $Q$. We say that a tiling of $\Si$ with congruent copies of $Q$ is
{\it trivial}, if every square $\si _{i,j}$ is tiled by four tiles of the tiling.
\begin{theorem}\label{t3}
Every tiling of $\Si$ with congruent copies of $Q$ is trivial.
\end{theorem}
We can deduce Lemma \ref{l6} from Theorem \ref{t3} as follows.
Suppose $Q$ is a reptile. Then there is a quadrilateral $Q'$ similar to $Q$
such that $Q'$ can be tiled with $k\ge 2$ congruent copies of $Q$.
Let the vertices of $Q'$ be $X,Y,V,W$ such that $\ol{XY}=\la a$, $\ol{YV}=\la$,
$\ol{VW}=\la$ and $\ol{WX}=\la b$, where $\la =\sqrt k >1$. We may assume that
$X$ is the origin, $Y$ is the point $(\la a,0)$ and $W$ is the point
$(0,\la b)$. The square $\si '_{0,0}=[0, \la (a+b)]\times  [0, \la (a+b)]$ can be
tiled with four copies of $Q'$, and thus by $4k$ congruent copies of $Q$.
Now $\Si$ is tiled with congruent copies of $\si '_{0,0}$. Translating the
tiling of $\si '_{0,0}$ into each of these squares we obtain a tiling of $\Si$
with congruent copies of $Q$.

By Theorem \ref{t3}, this tiling must be trivial. Let $i$ be the largest
integer with $i(a+b)< \la a$. Then the line $\ell$ going through the vertices
$Y$ and $V$ cuts the square $\si _{(i,0)}$ into two parts. Since the tiling is
trivial and $\ell \cap \si _{i,0}$ is covered by the boundaries of the tiles,
$\ell$ must intersect the $x$ axis at the point $i(a+b)+a$, and must
intersect the line $y=a+b$ at the point $(i(a+b)+b,a+b)$.
Then $\ell$ also cuts the square $\si _{i,1}$ into two parts.
However, $\ell \cap \si _{i,1}$ cannot be covered by the boundaries of the
tiles, since $b<a$, and thus $\ell$ enters the interior of a tile in
$\si _{i,1}$, which is impossible.

\begin{figure}
\includegraphics[width=4in]{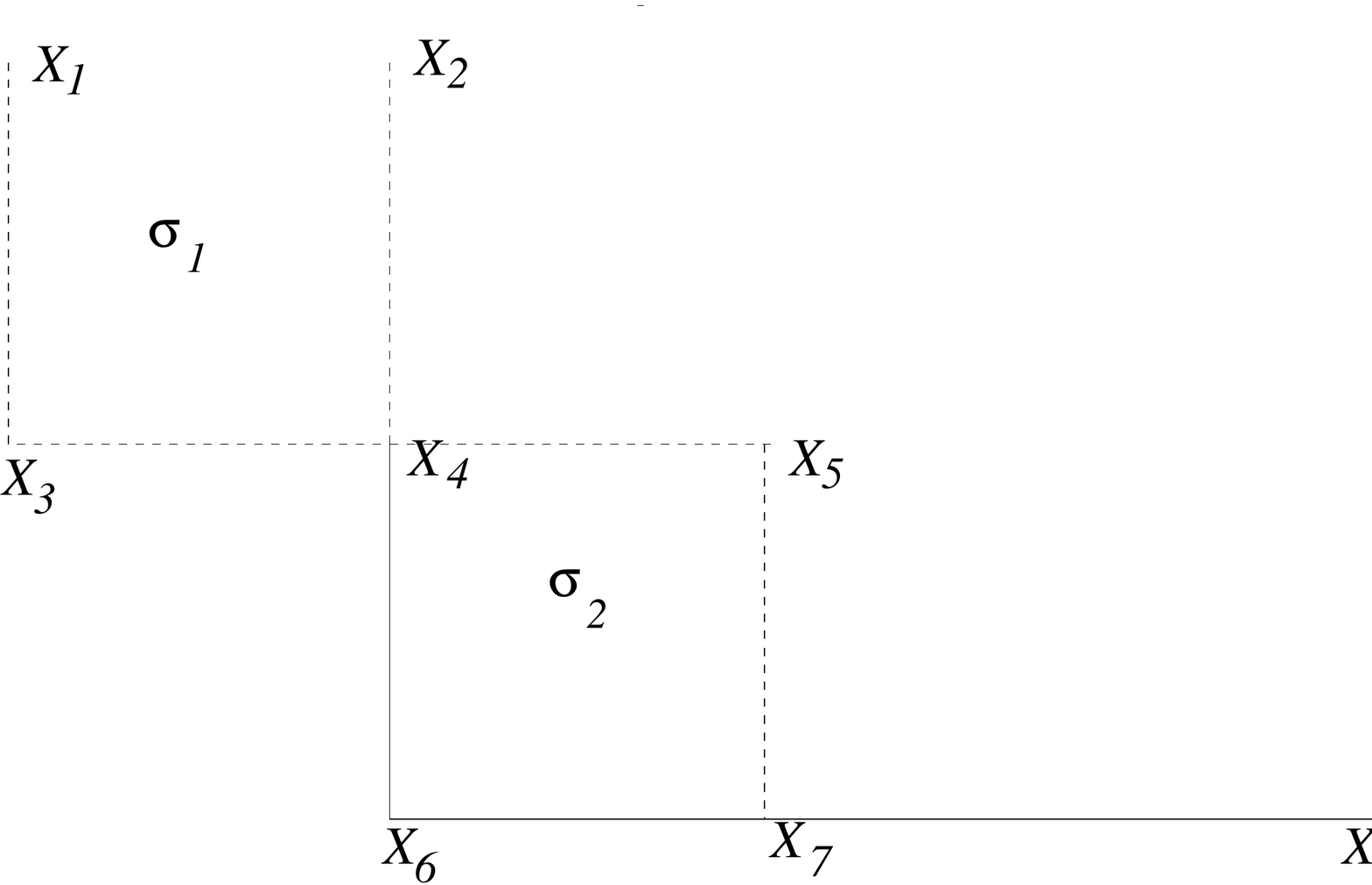}
\caption{}
\label{f8}
\end{figure}

The rest of the section is devoted to the proof of Theorem \ref{t3}.
Suppose a tiling of $\Si$ is given.
\begin{lemma}\label{l13}
Let $X_1 \stb X_8$ be points such that $X_3 ,X_4 ,X_2 ,X_1$ are the vertices of a
square $\si _1$  of side length $a+b$, $X_6 ,X_7 ,X_5 ,X_4$ are vertices of a
square $\si _2$  of side length $a+b$, and $\ol{X_7 X_8}>1$. (See Figure
\ref{f8}.) Suppose that at least one of the broken lines $X_1 X_3 X_4  X_6 X_8$
and $X_2 X_6 X_8$ is covered by $U$. Then $\si _2$ is tiled with four tiles.
\end{lemma}

\proof There is a tile $T_1 \su \si _2$ having a vertex at $X_6$. Let
the vertices of $T_1$ be $X_6 ,Y_1 ,Y_2 ,Y_3$ such that $Y_1$ is an inner
point of the segment $X_6 X_7$ and $Y_3$ is an inner point of the segment
$X_6 X_4$. By Lemma \ref{l12} we have $\ol{X_6 Y_1} =a$ and
$\ol{X_6 Y_3} =b$, or the other way around. We consider these two cases
separately.

\noi
{\bf Case I:} {\it $\ol{X_6 Y_1} =a$ and $\ol{X_6 Y_3} =b$.} Then there
is a tile $T_2$ with vertices $Y_3 Y_4 Y_5 Y_6$ such that $Y_4 \in \si _2$
and $Y_6$ is on the segment $X_6 X_2$. Since the angle of $T_2$ at $Y_3$
equals $\al$, we have $\ol{Y_3 Y_4} =a$ or $1$. We consider the two cases
separately. Let $\ell$ denote the line going through the points $Y_3$
and $Y_2$.

\begin{figure}
\includegraphics[width=5.5in]{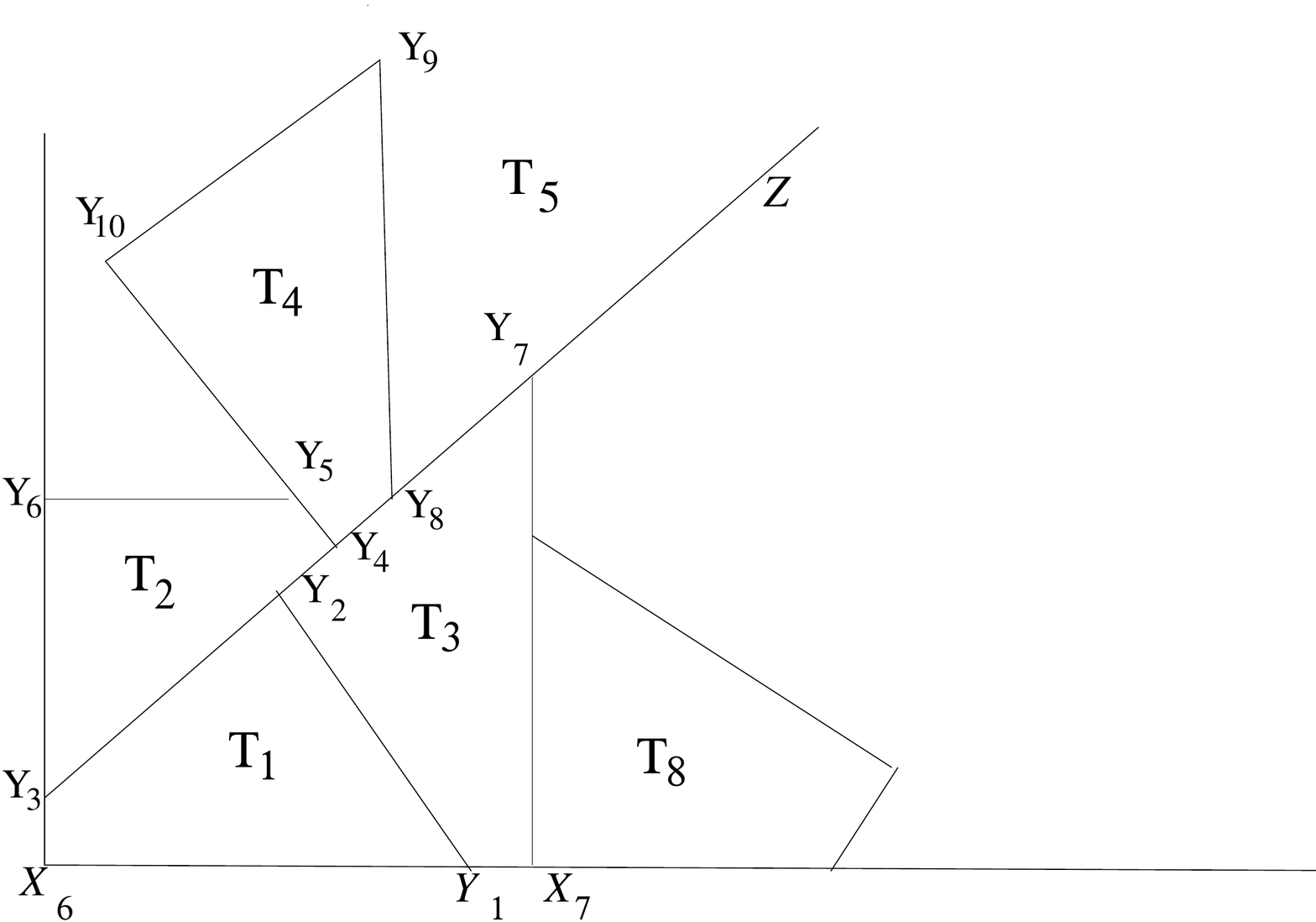}
\caption{}
\label{f9}
\end{figure}

\noi
{\bf Case I.1:} $\ol{Y_3 Y_4} =a$. Then $Y_4 Y_2 Y_1 X_7$ is a barrier, and
then, by Lemma \ref{l9}, there is a tile $T_3$ with vertices $Y_1 , X_7 ,Y_7 ,
Y_2$. (See Figure \ref{f9}.) We have $\ol{Y_2 Y_4}=a-1 <\ol{Y_2 Y_7}=1$, and
$\ol{Y_4 Y_7}=\ol{Y_3 Y_7}- \ol{Y_3 Y_4}=2-a$.

We show that there is a point $Z$ on $\ell$ such that $Y_3 Z \su U$ and $Y_7$
is an inner point of the segment $Y_3 Z$.

There is a tile $T_4$ with vertices $Y_4 ,Y_8, Y_{9}, Y_{10}$ such that $Y_8$
is on the line $\ell$. If $\ol{Y_4 Y_8} \ge 1$ then we take $Z=Y_8$.
On the other hand, if $\ol{Y_4 Y_8}=b$ then, by $2-a>b$, the point $Y_8$
is an inner point of the segment $Y_4 Y_7$. In this case, however, there is a 
tile $T_5$ with vertices $Y_8 ,Y_{11}, Y_{12}, Y_{13}$ such that $Y_{11}$
is on the line $\ell$. Since the angle of $T_5$ at $Y_8$ is $\al$, we have
$\ol{Y_8 Y_{11}}\ge 1$, and thus we can take $Z=Y_{11}$.

In both cases, $Z Y_7 X_7 X_8$ is a barrier as in Lemma \ref{l10}, so we
have a tile $T_8$ having sides of length $1$ on each of the segments $Y_7 X_7$
and $X_7 X_8$ (see Figure \ref{f9}). Since $\ol{Y_{7} X_{7}}=a>1$ and
$\ol{X_{7} X_{8}}>1$ by assumption, this contradicts Lemma \ref{l12}.

\noi
{\bf Case I.2:} $\ol{Y_3 Y_4} =1$. Then the vertices of $T_2$ are $Y_3 ,Y_2 ,
Y_5 ,X_4$. (That is, we have $Y_4 =Y_2$ and $Y_6 =X_4$.)

\begin{figure}
\includegraphics[width=5.5in]{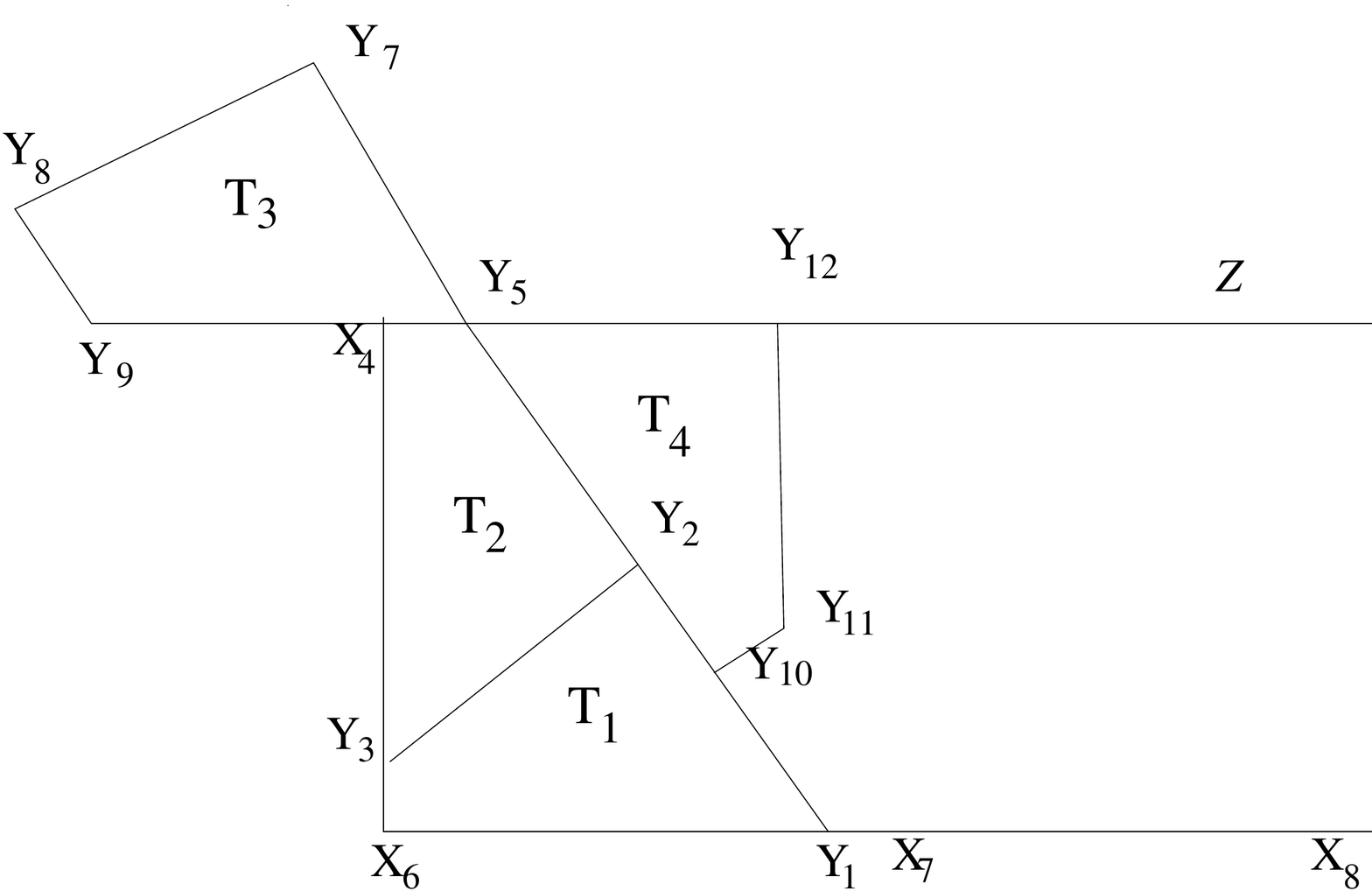}
\caption{}
\label{f10}
\end{figure}

We prove that there is a point $Z$ on the line going through $X_3 ,X_4$
and $Y_5$ such that the segment $X_4 Z$ is in $U$, and $Y_5$ is an inner
point of the segment $X_4 Z$. Suppose this is not true. Then there is
a tile $T_3$ with vertices $Y_5 ,Y_7 ,Y_8 ,Y_9$ as in Figure \ref{f10}.

Suppose that the angle of $T_3$ at $Y_5$ equals $\ga$. Since the angle of $T_2$
at $Y_5$ is also $\ga$, it follows that there are other tiles having a vertex
at $Y_5$. Since $2\pi -2\ga =2\al >\pi /2$, there must be two such tiles, each
of them having angle $\al$ at $Y_5$. From this it is clear that there is a
point $Z$ with the required properties; in fact, one of the vertices of the
tiles will do. Therefore, if there is no such point $Z$, then 
the angle of $T_2$ at $Y_5$ cannot be $\ga$.

If the broken line $X_2 X_6 X_8$ is in $U$, then we have $Y_9 =X_4$, and the
angle of $T_3$ at $Y_5$ is $\ga$, which is impossible.
Therefore, the broken line $X_1 X_3 X_4  X_6 X_8$ is in $U$, and thus
$X_1 X_3 Y_5 Y_7$ is a barrier having a right angle at $X_3$.
If the angle of $T_3$ at $Y_5$ is $\pi /2$, then, by (i) of Lemma \ref{l8},
$\ol{X_3 Y_5}$ equals a sum $c_1 +\ldots +c_n$, where each $c_i$ equals
one of $1+b, a+b, 2, 1+a$. Now $\ol{X_3 Y_5}=a+b+b$. Since $a+2b <2(1+b)$,
we have $n=1$, and thus one of the equalities $a+2b=1+b$, $a+2b=a+b$, $a+2b=2$,
$a+2b=1+a$ must hold. However, the cases $b=1/2$ and $a=2-2b$ are excluded,
so this is impossible.

Finally, suppose that the angle of $T_3$ at $Y_5$ is $\al$. Let
$S_1 \stb S_n$ be a tiling of the barrier $X_1 X_3 Y_5 Y_7$, and let
$X_3 =Z_0 \stb Y_n =Y_5$ be the corresponding partition. Then the angle of
$S_n$ at $Y_5$ is $\al$, and thus $\ol{Z_{n-1} Z_n}$ is either $1$ or $a$.
Therefore, $X_3 Z_{n-1}$ is the base of a barrier of two right angles.
Then $\ol{X_3 Z_{n-1}}$ equals either $a+2b-1<1+b$ or $2b<1+b$, contradicting
(i) of Lemma \ref{l8}.

We have proved the existence of a segment $X_3 Z \su U$ such that $Y_5$ is
an inner point of the segment $X_4 Z$. Then there is a tile $T_4$
with vertices $Y_5 ,Y_{10}, Y_{11}, Y_{12}$ such that $Y_{10}$ is on the segment
$Y_1 Y_5$. Since the angle of $T_4$ at $Y_5$ is $\al$, we have
$\ol{Y_5 Y_{10}}=a$ or $1$. If it equals $a$, then $Y_{11} Y_{10} Y_1 X_7$ is
a barrier with $\ol{Y_{10} Y_1}=2-a<1$. Thus the tiling of this barrier
has more than one tile, since $2-a$ is different from each of $a, b$ and $1$.
Thus $2-a\ge 1+b$, which is also impossible.

Therefore, we must have $\ol{Y_5 Y_{10}}=1$, and thus $Y_{10}=Y_2$ and
$Y_{12}=X_5$. Then $Y_{11} Y_2 Y_1 X_7$ is a barrier as in Lemma \ref{l9}.
Thus there is a tile with vertices $Y_2 ,Y_1 ,X_7 ,Y_{11}$ showing that
$\si _2$ is tiled with four tiles.

\noi
{\bf Case II:} {\it $\ol{X_6 Y_1} =b$ and $\ol{X_6 Y_3} =a$.} Then there
is a tile $T_2$ with vertices $Y_3 Y_4 Y_5 Y_6$ such that $Y_4$ is on the line
going through the points $Y_3$ and $Y_2$.
Since the angle of $T_2$ at $Y_3$ equals $\ga$, we have $\ol{Y_3 Y_4} =b$ or
$1$. We consider the two cases separately. 

\begin{figure}
\includegraphics[width=5in]{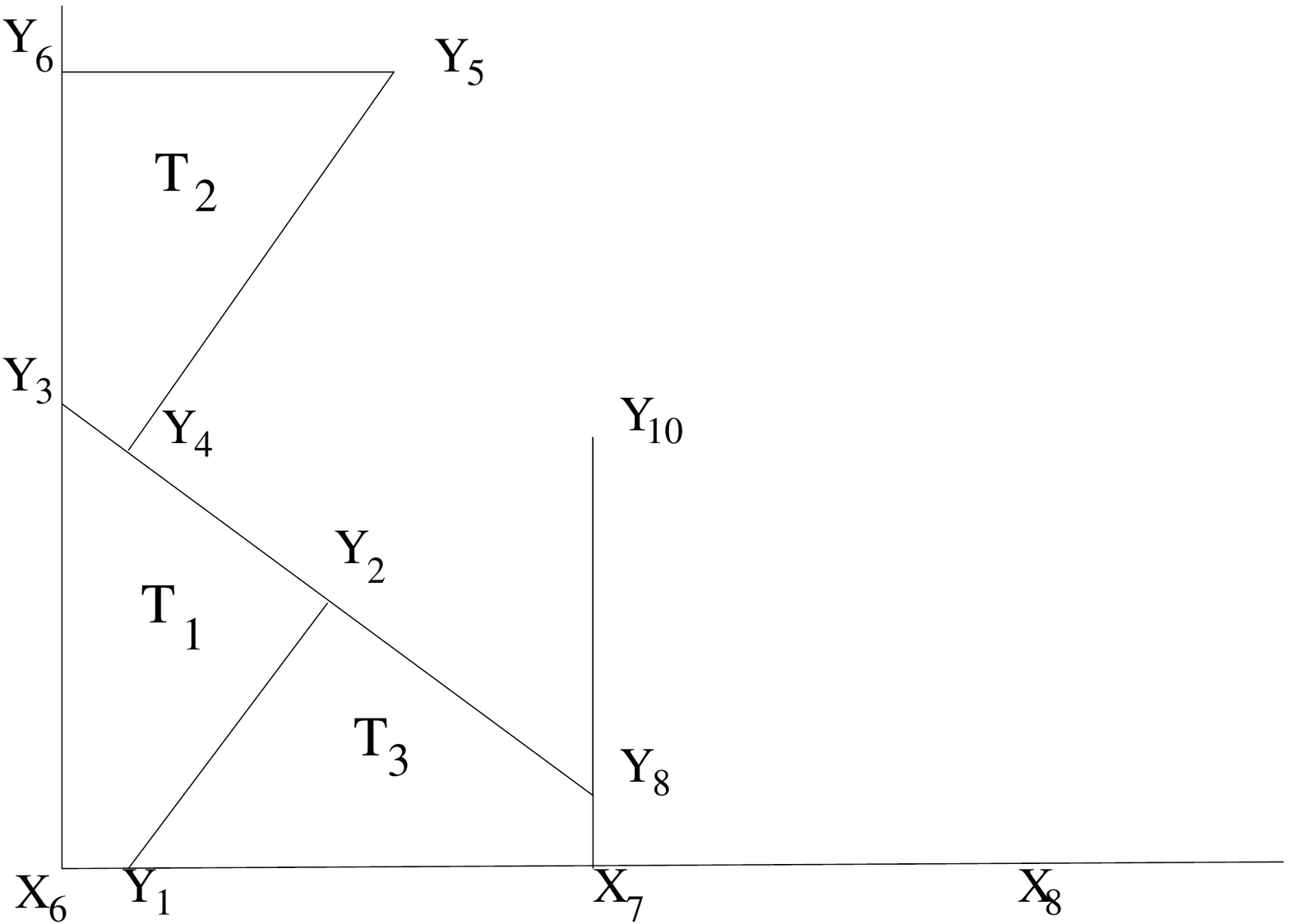}
\caption{}
\label{f11}
\end{figure}

\noi
{\bf Case II.1:} $\ol{Y_3 Y_4} =b$ (see Figure \ref{f11}). There is a tile
$T_3 \su \si _2$ with vertices $Y_1 ,Y_7 Y_8 Y_9$ such that $Y_7$ is on the
line $X_6 X_7$. Since the angle of $T_3$ at $Y_1$ equals $\al$, we have
$\ol{Y_1 Y_9} =a$ or $1$. If $\ol{Y_1 Y_9} =a$, then $Y_5 Y_4 Y_2 Y_9$ is a
barrier with two right angles and $\ol{Y_4 Y_2} =1-b<1$, which is impossible.
Therefore, we have $\ol{Y_1 Y_9} =1$, and then $Y_7 =X_7$ and $Y_9 =Y_2$.

By Lemma \ref{l11}, there exist a point $Y_{10}$ such that $Y_7 Y_{10} \su U$
and $Y_8$ is an inner point of the segment $Y_7 Y_{10}$. Then $Y_5 Y_4 Y_8 Y_{10}$
is a barrier with angles $\pi /2$ and $\al$, and with $\ol{Y_4 Y_8}=2-b$. Since
$2-b$ is different from each of $b,1$ and $a$, the tiling of the barrier
consists of more than one tile. One of the tiles has angle $\al$ at $Y_8$, so
its side on the base has length $a$ or $1$. The other tiles tile a barrier with
two right angles, and thus $2-b=\ol{Y_4 Y_8} \ge (1+b)+1$, which is absurd. So
this case cannot happen.

\noi
{\bf Case II.2:} $\ol{Y_3 Y_4} =1$. Then $Y_4 =Y_2$ and $Y_6 =X_4$.
There is a tile $T_3 \su \si _2$ with vertices $Y_1 ,Y_7 Y_8 Y_9$ such that
$Y_7$ is on the line $X_6 X_7$. Since the angle of $T_3$ at $Y_1$ equals $\al$,
we have $\ol{Y_1 Y_7} =a$ or $1$. Let $\ell$ denote the line going through
the points $Y_1 , Y_2 , Y_5$.

\begin{figure}
\includegraphics[width=3in]{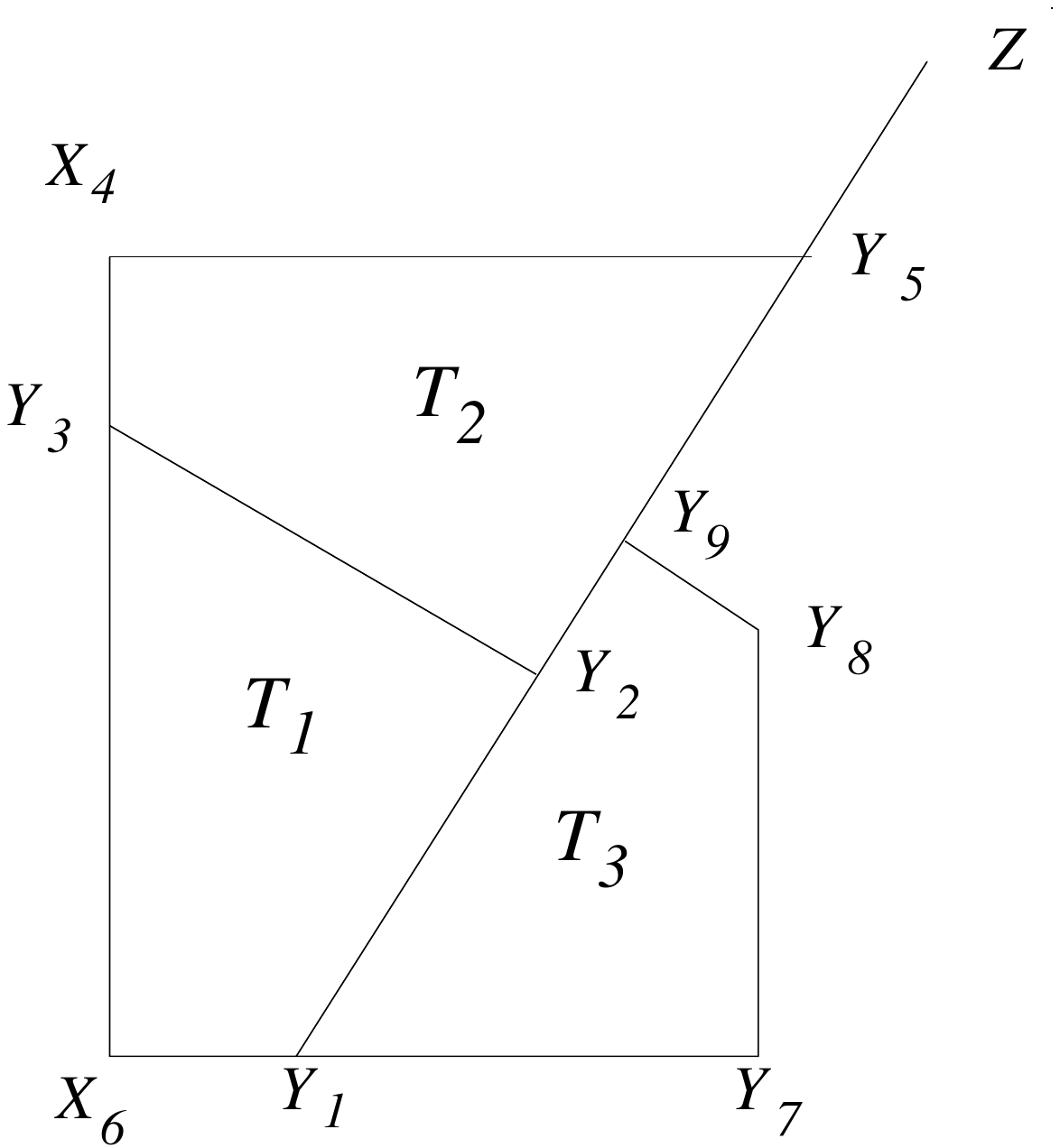}
\caption{}
\label{f13}
\end{figure}

Suppose $\ol{Y_1 Y_7} =1$ (see Figure \ref{f13}). We prove that
there exists a point $Z$ on the line $\ell$ such that $Y_1 Z \su U$ and $Y_5$
is an inner point of the segment $Y_1 Z$.

There is a tile $T_4$ having a vertex at $Y_9$. Let the vertices of $T_4$ be
$Y_9 , Y_{10} , Y_{11} , Y_{12}$, where
$Y_{12}$ is on the line $\ell$. If $\ol{Y_2 Y_{12}}\ge 1$, then we take
$Z=Y_{12}$. If $\ol{Y_2 Y_{12}}=b$, then $T_4$ has angle $\ga$ at $Y_{12}$,
and then there is a tile $T_5$ having angles $\al$ at $Y_{12}$. In this case
let $Z$ be the vertex of $T_5$ lying on $\ell$.

If the broken line $X_2 X_6 X_8$ is covered by $U$ then, by Lemma \ref{l10},
there is a tile $T_4$ such that $T_4$ has a vertex at $X_4$, and has sides of
length $1$ on both segments $X_3 X_4$ and $X_4 X_5$. However,
these segments are longer than $1$, which contradicts Lemma \ref{l12}.

If the the broken line $X_1 X_3 X_4  X_6 X_8$ is covered by $U$, then 
$X_1 X_3 Y_5 Z$ is a barrier such that $X_1 X_3 Y_5 \szog =\pi /2$,
$X_3 Y_5 Z \szog =\ga$ and $\ol{X_3 Y_5}=2a+b$. By lemma \ref{l9a}, this
implies $a+b=\ol{X_1 X_3}\le 1$, which is not true.

Therefore, this case cannot happen, and we have $\ol{Y_1 Y_7} =a$.
Applying Lemma \ref{l11} we find a point $W$ on the line going
through the points $X_7 , Y_8$ and $X_5$ and such that $X_7 W\su U$ and
$Y_8$ is an inner point of the segment $X_7 W$.
Then $Y_5 Y_9 Y_8 W$ is a barrier such that $Y_5 Y_9 Y_8 \szog =\pi /2$,
$Y_9 Y_8 W \szog =\al$, and $\ol{Y_2 Y_8}=1$. Clearly, $Y_5 Y_9 Y_8 W$ is
tiled with a single tile $T_4$ such that $T_1 ,T_2 ,T_3$
and $T_4$ tile $\si _2$. \hfill $\square$

Now we turn to the last step of the proof of Theorem \ref{t3}. Let $X_{i,j}$
denote the point with coordinates $(i(a+b),j(a+b))$. We denote by
$\Xi$ the set of squares $\si _{i,j}$ that are tiled with four tiles.
We have to prove that $\si _{i,j} \in \Xi$ for every $i,j\ge 0$.

We prove this statement by induction on $k=i+j$. By Lemma \ref{l13},
$\si _{0,0}\in \Xi$, since the segments $X_{0,0} X_{0,2}$ and $X_{0,0} X_{2,0}$
belong to $U$.

We have $\si _{1,0}\in \Xi$, since the broken line $X_{0,2} X_{0,1} X_{1,1}
X_{1,0} X_{3,0}$ belongs to $U$, and Lemma \ref{l13} applies.

We also have $\si _{0,1}\in \Xi$, since the broken line $X_{0,3} X_{0,1} X_{2,1}$
belongs to $U$, and Lemma \ref{l13} applies.

Let $k\ge 2$, and suppose that $\si _{i,j} \in \Xi$ for every $i,j$
such that $i+j<k$.

We prove $\si _{k-j,j} \in \Xi$ by induction on $j$. As for $j=0$,
$\si _{k,0} \in \Xi$ follows from Lemma \ref{l13}, since $\si _{k-1,0} ,
\si _{k-2,1} \in \Xi$ by the induction hypotheses, and thus the broken line
$X_{k-1,2} X_{k-1,1} X_{k,1} X_{k,0}X_{k+2,0}$ is covered by $U$.

Let $0<j\le k$, and suppose that $\si _{k-j+1,j-1} \in \Xi$. Then there
are three cases.

If $j\le k-2$, then $\si _{k-j-2,j+1}$, $\si _{k-j-1,j}$, $\si _{k-j, j-1}$,
$\si _{k-j+1,j-1} \in \Xi$ by the induction hypotheses. Therefore, the broken
line
$$X_{k-j-1,j+2} X_{k-j-1,j+1} X_{k-j,j+1} X_{k-j,j} X_{k-j+2,j}$$
is covered by $U$, and Lemma \ref{l13} applies.

If $j=k-1$, then $\si _{0,k-1} ,\si _{1,k-2}$, $\si _{2,k-2} \in \Xi$
by the induction hypotheses. Therefore, the broken line $X_{0,k+1}
X_{0,k} X_{1,k} X_{1,k-1} X_{3,k-1}$ is covered by $U$, and Lemma \ref{l13}
applies.

Finally, If $j=k$, then $\si _{1,k-1}  \in \Xi$ by the induction hypotheses.
In this case the broken line $X_{0,k+2} X_{0,k} X_{2,k} $ is covered by $U$, and
Lemma \ref{l13} applies.

This completes the proof of Theorem \ref{t3} and that
of Theorem \ref{t1}. \hfill $\square$

\vfill \eject

\end{document}